\documentclass[12pt,twoside,leqno,openany]{amsart}
\usepackage{amssymb,amsbsy,amsmath,amsfonts,amscd,times}
\usepackage{graphics,color,xy,xypic,footmisc,fancyhdr,multicol}
\usepackage{fancybox,graphicx,mathrsfs}
\usepackage[T1]{fontenc}
\usepackage{breqn}
\let\mathcal\mathscr
\newtheorem{theorem}{Theorem}
\newtheorem{lemma}{Lemma}

\newcommand{\C}{\mathbb{C}}

\newcommand{\R}{\mathbb{R}}

\newcommand{\ov}{\overline}
\newcommand{\LL}{\mathcal{L}_1}
\newcommand{\Lb}{\overline{\mathcal{L}_1}}
\newcommand{\T}{\mathcal{T}}
\newcommand{\s}{\mathcal{S}}

\newcommand{\zb}{\ov{z}}
\newcommand{\A}{{\sf a}}
\newcommand{\ab}{\ov{{\sf a}}}
\newcommand{\bb}{{\sf b}}
\newcommand{\bbb}{\ov{{\sf b}}}
\newcommand{\cc}{{\sf c}}
\newcommand{\cb}{\ov{{\sf c}}}
\newcommand{\dd}{{\sf d}}
\newcommand{\db}{\ov{{\sf d}}}
\newcommand{\ee}{{\sf e}}
\newcommand{\eb}{\ov{{\sf e}}}
\newcommand{\ff}{{\sf f}}
\newcommand{\fb}{\ov{{\sf f}}}
\newcommand{\G}{{\sf g}}
\newcommand{\hh}{{\sf h}}

\newcommand{\kk}{{\sf k}}

\newcommand{\LC}{\sf LC}

\newcommand{\KK}{\mathcal{K}}
\newcommand{\Kb}{\ov{\mathcal{K}}}

\newcommand{\smallbullet}{{\scriptscriptstyle{\bullet}}}

\begin{document}

\setcounter{page}{116}

\title[]{
Lie algebras of infinitesimal automorphisms for the model manifolds of general classes ${\sf II}$, ${\sf III_2}$ and ${\sf IV_2}$}

\author{Samuel Pocchiola}
\address{Samuel Pocchiola ---  D\'epartement de math\'ematiques d'Orsay, b\^atiment 425, Facult\'e des sciences d'Orsay,
Universit\'e Paris-Sud, F-91405 Orsay Cedex, france}
\email{samuel.pocchiola@math.u-psud.fr}

\maketitle

\bigskip
\section*{abstract}
We determine the Lie algebra of infinitesimal CR-automorphisms of the 
models of general classes ${\sf II}$,  ${\sf III_2}$ and ${\sf IV_2}$ through Cartan's equivalence method. 

\section{Introduction}

The classification of CR-manifolds up to dimension $5$ has highlighted the existence of $6$ non-trivial classes of CR-manifolds,
which have been referred to as general classes ${\sf I}$, ${\sf II}$, ${\sf III}_1$, ${\sf III}_2$, ${\sf IV}_1$ and ${\sf IV}_2$ \cite{MPS}.
Each of these classes entails a distinguished manifold, the model, whose Lie algebra of infinitesimal CR-automorphisms is of maximal dimension.  
It plays a special role, as CR-manifolds belonging to the same class can be viewed as its deformations, generally by the way of Cartan connection.
The aim of this paper is to determine the Lie algebra of infinitesimal CR-automorphisms
of the models for general classes ${\sf II}$,  ${\sf III}_2$ and ${\sf IV}_2$. This is already known \cite{BES-2007, Kaup-Zaitsev}
for general classes ${\sf II}$ (Engel manifolds) 
and ${\sf IV}_2$ (2-nondegenerate, $5$-dimensional CR-manifolds of constant Levi rank $1$), but is unknown, to our knowledge, in the case 
of general class 
${\sf III}_2$. In our view, the main interest of this paper is to provide a unified treatment for the $3$ classes through
the use of Cartan's equivalence method, in the spirit of \cite{Olver-1995}.
Cartan's equivalence method has indeed been employed recently
to solve the 
equivalence problem for general classes ${\sf II}$, ${\sf III}_2$ and ${\sf IV_2}$ \cite{pocchiola, pocchiola2, pocchiola3}. For each of these classes,
the solution to the equivalence problem for the model has been of a great help for the treatment of the general case,
as a similar structure of normalizations of the group parameters occurs in both cases.

For general class ${\sf II}$, the model is provided by Beloshapka's cubic in $\C^3$, which is the CR-manifold defined by the equations:
\begin{equation*}
{\sf B}: \qquad \qquad
\begin{aligned}
w_1 & = \ov{w_1} + 2 \, i \, z \ov{z}, \\
w_2 & = \ov{w_2} + 2 \, i \, z \ov{z} \left( z + \ov{z} \right).
\end{aligned}
\end{equation*}
For general class ${\sf III}_2$, the model is the $5$-dimensional submanifold ${\sf N} \subset \C^4$ defined by:
\begin{equation*}
{\sf N}: \qquad \qquad
\begin{aligned}
w_1 & = \ov{w_1} + 2 \, i \, z \ov{z}, \\
w_2 & = \ov{w_2} + 2 \, i \, z \ov{z} \left( z + \ov{z} \right), \\
w_3 & = \ov{w_3} + 2 i \, z \zb ( z^2 + \frac{3}{2} \, z \zb + \zb^2 ).
\end{aligned}
\end{equation*}
For general class ${\sf IV}_2$, the model is provided by the tube over the future light cone, ${\sf LC} \subset \C^3$, defined by:
\begin{equation*}
{\sf LC}: \,\,\,\,\,\,\,\,\,\,\,\,\,\,\,\,\,\,\,\,\,\,\,\,\, 
\left( {\sf Re} \,z_1 \right)^2 -\left( {\sf Re} \,z_2 \right)^2 
- \left( {\sf Re} \,z_3 \right)^2
=
0, \qquad \qquad {\sf Re} \,z_1 > 0
.\end{equation*}

A Cartan connection has been constructed for CR-manifolds belonging to general class ${\sf II}$ \cite{BES-2007, pocchiola2} and
${\sf III}_2$ \cite{pocchiola3}. The equivalence problem for manifolds belonging to general class ${\sf IV}_2$ has been solved either by the determination of
an absolute parallelism 
\cite{ Isaev-Zaitsev, pocchiola}, or the construction of a Cartan connection \cite{Medori-Spiro}.
We use Cartan's equivalence method for which we refer to \cite{Olver-1995} 
as a standard reference.

\section{Class ${\sf II}$}
This section is devoted to the determination of the Lie algebra of CR-automorphisms
of Beloshapka's cubic in $\C^3$, which is the CR-manifold defined by the equations:
\begin{equation*}
{\sf B}: \qquad \qquad
\begin{aligned}
w_1 & = \ov{w_1} + 2 \, i \, z \ov{z}, \\
w_2 & = \ov{w_2} + 2 \, i \, z \ov{z} \left( z + \ov{z} \right).
\end{aligned}
\end{equation*}
It is the model manifold 
for generic $4$-dimensional CR-manifolds of CR 
dimension $1$ and real codimension $2$, i.e. CR-manifolds belonging to class ${\sf II}$, 
in the sense that any such manifold might be viewed
as a deformation of Beloshapka's cubic by the way of a Cartan connection \cite{BES-2007, pocchiola2}.
The main result of this section is:

\begin{theorem}
\label{thm:Bc}
Beloshapka's cubic,
\begin{equation*}
{\sf B}: \qquad \qquad
\begin{aligned}
w_1 & = \ov{w_1} + 2 \, i \, z \ov{z}, \\
w_2 & = \ov{w_2} + 2 \, i \, z \ov{z} \left( z + \ov{z} \right),
\end{aligned}
\end{equation*}
has a ${\bf 5}$-dimensional Lie algebra of CR-automorphisms. 
A basis for the Maurer-Cartan forms of ${\sf aut_{CR}}({\sf B})$ is provided
by the $5$ differential $1$-forms  $\sigma$, $\rho$, $\zeta$,  $\ov{\zeta}$, $\alpha$,
which satisfy the structure equations:
\begin{equation*}
\begin{aligned}
d \sigma &= 3 \left. \alpha \wedge \sigma \right. + \left. \rho \wedge \zeta \right. 
+ \left.  \rho \wedge \ov{\zeta} \right., \\
d \rho &  =  2 \left. \alpha \wedge \rho \right. + i \, \left. \zeta \wedge \ov{\zeta} \right., \\
d \zeta &= \left. \alpha \wedge \zeta \right., \\
d \ov{\zeta} &= \left. \alpha \wedge \ov{\zeta} \right., \\
d \alpha & = 0.
\end{aligned}
\end{equation*}
\end{theorem}

\subsection{Initial G-structure.}
 
The vectors field $\LL$ defined by:
\begin{equation*}
\LL:=
\frac{\partial}{\partial z}
+
i \, \zb \,
\frac{\partial}{\partial u_1}
+ i \,
\left( 2 z \zb + \zb^2 \right) \, \frac{\partial}{\partial u_2},
\end{equation*} 
together with its conjugate:

\begin{equation*}
\Lb:= \frac{\partial}{\partial \zb}
- i \, z \,
\frac{\partial}{\partial u_1}
- i \,
\left( 2 z \zb + z^2 \right) \, \frac{\partial}{\partial u_2},
\end{equation*}
constitute a basis of $T^{1,0}_p {\sf B}$ at each point $p$ of ${\sf B}$.
Moreover the vector fields $\T$ and $\s$ defined by:
\begin{equation*}
\T: = i \, \big[\LL, \Lb] 
,\end{equation*}
and
\begin{equation*}
\s: = \big[\LL, \T \big], 
\end{equation*}
complete a frame on ${\sf B}$:
\begin{equation*}
\big\{\mathcal{S},\,
\mathcal{T},\,\mathcal{L}, \,
\overline{\mathcal{L}}
\big\}
.
\end{equation*}
The expressions of $\T$ and $\s$ are:

\begin{equation*}
\begin{aligned}
\T &  = 2 \, \frac{\partial}{\partial u_1} + \left( 4z + 4 \zb \right) \frac{\partial}{\partial u_2}, \\
\s &  = 4 \, \frac{\partial}{\partial u_2}
.\end{aligned}
\end{equation*}
The dual coframe $ \big(\sigma_0,\,
\rho_0,
\,
\zeta_0,
\,
\overline{\zeta_0}
\big)$
is thus given by:

\begin{equation*}
\begin{aligned}
\sigma_0 &= \frac{i}{4} \, \zb^2 \, dz - 
\frac{i}{4} \, z^2 \, d \zb -  \left( \frac{1}{2} \, z + \frac{1}{2} \, \zb \right) du_1 + \frac{1}{4} \, du_2 
,\\
\rho_0 &=  - \frac{i}{2} \, \zb \,  dz + \frac{i}{2} \, z \, d\zb + \frac{1}{2} \, du_1, \\
\zeta_0 &= dz,\\
\ov{\zeta_0} &= d\zb.
\end{aligned}
\end{equation*}
We deduce the structure equations enjoyed by
$
\big(\sigma_0,\,
\rho_0,\,
\zeta_0, \,
\overline{\zeta_0},
\big)
$:
\begin{equation}
\label{eq:B}
\begin{aligned}
d \sigma_0 & = \rho_0 \wedge \zeta_0 + \rho_0 \wedge \ov{\zeta_0}, \\
d \rho_{0} & =  i \, \zeta_0 \wedge \ov{\zeta_0},\\
d \zeta_0 & = 0, \\
d \ov{ \zeta_{0}} & = 0.
\end{aligned}
\end{equation}

As the torsion coefficients of these structure equations are constants, we have the following result:

\begin{lemma}
Beloshapka's cubic
is locally isomorphic to a Lie group whose Maurer-Cartan forms satisfy the structure equations 
$(\ref{eq:B})$.
\end{lemma} 

The matrix Lie group which encodes suitably the equivalence problem for Beloshapka's cubic (see \cite{pocchiola2}) is the 
$10$-dimensional Lie group $G_1$ whose elements $g$ are of the form:
\begin{equation*}
g := \begin{pmatrix}
{\A^2} \ab & 0 & 0 & 0  \\
\cc & \A \ab & 0 & 0 \\
\dd & \bb & \A & 0 \\
\ee & \bbb & 0 & \ab
\end{pmatrix}.
\end{equation*}

With the notations:
\begin{alignat*}{1}
\omega_{0} :=
\begin{pmatrix}
\sigma_{0} \\ \rho_0 \\ \zeta_{0} \\ \ov{\zeta}_0
\end{pmatrix}, & \qquad
\omega:=
\begin{pmatrix}
\sigma \\ \rho \\ \zeta \\ \ov{\zeta}
\end{pmatrix},
\end{alignat*}
we introduce the $G_1$-structure $P^1$ on ${\sf B}$ constituted 
by the coframes $\omega$ which satisfy the relation:

\begin{equation*}
\omega:= g \cdot \omega_0.
\end{equation*}
The proof of theorem (\ref{thm:Bc}) relies on successive reductions of $P^1$ through Cartan's equivalence method.

\subsection{Normalization of {\A}}
The structure equations for the lifted coframe $\omega$ are related to those of the base coframe $\omega_0$
 by the relation:
\begin{equation} 
\label{eq:str}
d \omega =  dg \cdot g^{-1} \wedge \omega + g \cdot d \omega_{0}.
\end{equation}
The term $ dg \cdot g^{-1} \wedge \omega$ 
depends only on the structure equations of $G_1$ and is expressed through its Maurer-Cartan forms.
The term $ g \cdot d \omega_{0}$  
contains the so-called torsion coefficients of the $G_1$-structure. We can compute it easily in terms of the forms 
$\sigma$, $\rho$,  $\zeta$, $\ov{\zeta}$,
 by a simple multiplication by $g$ in the formulae $(\ref{eq:B})$ and a linear change of variables.
The Maurer-Cartan forms for the group $G_1$ 
are given by the linearly independent entries of the matrix $dg \cdot g^{-1}$, which are:
\begin{equation*}
\begin{aligned}
\alpha^1  &: = {\frac {{ d\A}}{\A}}, \\
\alpha^2 &:= -{\frac {\bb{ d\A}}{{\A}^{2}{ \ab}}}+{\frac {{ d\bb}}{\A{ \ab}}}, \\
\alpha^3  &:= -{\frac {\cc{ d\A}}{{ \ab}\,{\A}^{3}}}-{\frac {\cc{ d\ab
}}{{{ \ab}}^{2}{\A}^{2}}} + { \frac {{ d\cc}}{{\A}^{2}{ \ab}}}, \\
\alpha^4 &: =-{\frac { \left( \dd\A{ \ab}-\bb\cc \right) { d\A}}{{\A}^{4}{
{ \ab}}^{2}}}-{\frac {\cc{ d\bb}}{{\A}^{3}{{ \ab}}^{2}}}+{\frac {{
 d\dd}}{{\A}^{2}{ \ab}}}, \\
\alpha^5 & :=-{\frac { \left( \ee\A{ \ab}-{ \bbb}\,\cc \right) { d\ab}
}{{\A}^{3}{{ \ab}}^{3}}}-{\frac {\cc{ d\bbb}}{{\A}^{3}{{ \ab}}^{2}}}+
{\frac {{ d\ee}}{{\A}^{2}{ \ab}}}, 
\end{aligned}
\end{equation*}
together with their conjugates.

The first structure equation is given by:

\begin{equation*}
d \sigma = 2 \left.\alpha^1 \wedge \sigma \right. + \left.\ov{\alpha^1} \wedge \sigma \right.
+ \left( \frac{\ee}{\A \aa^2} + \frac{\dd}{\A^2 \ab} \right) \left.\sigma \wedge \rho \right.
- \frac{\cc}{\A^2 \ab} \, \left.\sigma \wedge \zeta \right. 
-  \frac{\cc}{\A \ab^2} \, \left. \sigma \wedge \ov{\zeta} \right. + \left.\rho \wedge \zeta \right. + 
\frac{\A}{\ab} \, \left. \rho \wedge \ov{\zeta} \right..
\end{equation*}
from which we immediately deduce that $\frac{\A}{\ab}$ is an essential torsion coefficient which might
 be normalised to $1$ by setting:
\begin{equation*}
\A = \ab.
\end{equation*}

\subsection{Normalizations of $\bb$ and $\cc$}

We have thus reduced the $G_{1}$ equivalence problem on ${\sf B}$ 
to a $G_2$ equivalence problem, where $G_2$ is the $9$ dimensional real matrix 
Lie group whose elements are of the form
\begin{equation*}
g := \begin{pmatrix}
{\A}^3 & 0 & 0 & 0  \\
\cc & \A^2 & 0 & 0 \\
\dd & \bb & \A & 0 \\
\ee & \bbb & 0 & \A
\end{pmatrix},
\qquad \qquad \A \in \R.
\end{equation*}
The Maurer-Cartan forms of $G_2$ are given by:
\begin{equation*}
\begin{aligned}
\beta^1  &: = {\frac {{ d\A}}{\A}}, \\
\beta^2 &:= -{\frac {\bb d\A }{{\A}^{3}}}+{\frac {{ d\bb}}{\A^2}}, \\
\beta^3  &:=- 2 \, {\frac {\cc  d\A }{{\A}^{4}}} + { \frac {{ d\cc}}{{\A}^{3}}}, \\
\beta^4 & :=-{\frac { \left( \dd\A^2 -\bb \cc \right) { d\A}}{{\A}^{6}}}-{\frac {\cc{ d\bb}}{{\A}^{5}}}+{\frac {{
 d\dd}}{{\A}^{3}}}, \\
\beta^5 & :=-\frac { \left( \ee\A^2-{ \bbb}\,\cc \right)  d\A
}{\A^6}-\frac {\cc{ d\bbb}}{\A^{5}}+
\frac { d\ee}{\A^{3}}, \\
\end{aligned}
\end{equation*}
together with $\ov{\beta^2}$, $\ov{\beta^3}$, $\ov{\beta^4}$, $\ov{\beta^5}$.
Using formula $(\ref{eq:str})$, we get the structure 
equations for the lifted coframe $(\sigma, \rho, \zeta, \ov{\zeta})$ from those of the base coframe
$(\sigma_0, \rho_0, \zeta_0, \ov{\zeta_0})$ by a
 matrix multiplication and a linear change of coordinates, as in the first step:

\begin{multline*}
d \sigma = 3 \, \beta^1 \wedge \sigma  \\
+
U^{\sigma}_{\sigma \rho} \left. \sigma \wedge \rho \right.
+
U^{\sigma}_{\sigma \zeta} \left. \sigma \wedge \zeta \right.
+
U^{\sigma}_{\sigma \ov{\zeta}} \left. \sigma \wedge \ov{\zeta} \right.
+
\rho \wedge \zeta
+
\rho \wedge \ov{\zeta},
\end{multline*}

\begin{multline*}
d \rho
=
2 \beta^{1} \wedge \rho + \beta^3 \wedge \sigma \\
+
U^{\rho}_{\sigma \rho} \, \sigma \wedge \rho
+
U^{\rho}_{\sigma \zeta} \, \sigma \wedge \zeta
+
U^{\rho}_{\sigma \ov{\zeta}} \, \sigma \wedge \ov{\zeta}  \\
+
U^{\rho}_{\rho \zeta} \, \rho \wedge \zeta  
+
U^{\rho}_{\rho \ov{\zeta}} \, \rho \wedge \ov{\zeta} 
+
i \, \zeta \wedge \ov{\zeta}
,\end{multline*}

\begin{multline*}
d \zeta
=
{\beta}^{1} \wedge \zeta + {\beta}^{2} \wedge \rho + {\beta}^{4} \wedge \sigma \\
+
U^{\zeta}_{\sigma \rho} \, \sigma \wedge \rho
+
U^{\zeta}_{\sigma \zeta} \, \sigma \wedge \zeta 
+
U^{\zeta}_{\sigma \ov{\zeta}} \, \sigma \wedge \ov{\zeta} 
+
U^{\zeta}_{\rho \zeta} \, \rho \wedge \zeta \\
+
U^{\zeta}_{\rho \overline{\zeta}} \, \rho \wedge \overline{\zeta}
+
U^{\zeta}_{\zeta \overline{\zeta}} \, \zeta \wedge \overline{\zeta}
.\end{multline*}

We now proceed with the absorption phase. We introduce the modified Maurer-Cartan forms:
\begin{equation*}
\widetilde{\beta}^i= \beta^i  - y_{\sigma} \, \sigma - y_{\rho}^i \, \rho  - y_{\zeta}^i \, \zeta \, 
- y_{\overline{\zeta}}^i \, \overline{\zeta},
\end{equation*}
such that the structure equations rewrite:

\begin{multline*}
d \sigma = 3 \left. \widetilde{\beta}^1 \wedge \sigma \right. \\
+
\left( U^{\sigma}_{\sigma \rho} - 3 \, y^1_{\rho} \right) \, \left. \sigma \wedge \rho \right. 
+
\left( U^{\sigma}_{\sigma \zeta} -3 \,y^1_{\zeta} \right) \, \left. \sigma \wedge \zeta \right. \\
+
\left( U^{\sigma}_{\sigma \ov{\zeta}}  -3 \,y^1_{\ov{\zeta}} \right) \, \left. \sigma \wedge \ov{\zeta} \right. 
\rho \wedge \zeta
+
\rho \wedge \ov{\zeta},
\end{multline*}

\begin{multline*}
d \rho
=
2 \widetilde{\beta}^{1} \wedge \rho + \widetilde{\beta}^3 \wedge \sigma \\
+
\left( U^{\rho}_{\sigma \rho} + 2 \, y^1_{\sigma} - y^3_{\rho} \right) \, \left. \sigma \wedge \rho \right.
+
\left( U^{\rho}_{\sigma \zeta} - y^3_{\zeta} \right) \, \left. \sigma \wedge \zeta \right. \\
+
\left( U^{\rho}_{\sigma \ov{\zeta}} - y^3_{\ov{\zeta}} \right) \, \left. \sigma \wedge \ov{\zeta} \right.
+
\left( U^{\rho}_{\rho \zeta} - 2 \, y^1_{\zeta} \right) \, \left. \rho \wedge \zeta \right.  \\ 
+
\left( U^{\rho}_{\rho \ov{\zeta}} - 2 \, y^1_{\ov{\zeta}} \right) \, \left. \rho \wedge \ov{\zeta} \right.
+
i \, \left. \zeta \wedge \ov{\zeta} \right.
,\end{multline*}

\begin{multline*}
d \zeta
=
\widetilde{\beta^{1}} \wedge \zeta + \widetilde{\beta^{2}} \wedge \rho + \widetilde{\beta^{4}} \wedge \sigma \\
+
\left( U^{\zeta}_{\sigma \rho} + y^2_{\sigma}- y^4_{\rho} \right) \, \left. \sigma \wedge \rho \right.
+
\left( U^{\zeta}_{\sigma \zeta} + y^1_{\sigma} - y^4_{\zeta} \right) \, \left. \sigma \wedge \zeta \right. \\ 
+
\left( U^{\zeta}_{\sigma \ov{\zeta}} - y^4_{\ov{\zeta}} \right) \, \left. \sigma \wedge \ov{\zeta} \right. 
+
\left( U^{\zeta}_{\rho \zeta} + y^1_{\rho} - y^2_{\zeta} \right) \, \left. \rho \wedge \zeta \right. \\
+
\left( U^{\zeta}_{\rho \overline{\zeta}} - y^2_{\ov{\zeta}} \right) \, \left. \rho \wedge \overline{\zeta} \right.
+
\left( U^{\zeta}_{\zeta \overline{\zeta}} - y^1_{\ov{\zeta}} \right) \, \left. \zeta \wedge \overline{\zeta} \right.
.\end{multline*}

We get the following absorbtion equations:

\begin{alignat*}{3}
3 \, y^1_{\rho} & =  U^{\sigma}_{\sigma \rho}, & \qquad \qquad  
3 \,y^1_{\zeta} & =  U^{\sigma}_{\sigma \zeta}, & \qquad \qquad
3 \,y^1_{\ov{\zeta}} & =  U^{\sigma}_{\sigma \ov{\zeta}}, \\
-2 \, y^1_{\sigma} + y^3_{\rho}  &= U^{\rho}_{\sigma \rho}, &\qquad \qquad 
y^3_{\zeta}  &=  U^{\rho}_{\sigma \zeta}  , &\qquad \qquad   
 y^3_{\ov{\zeta}}   &= U^{\rho}_{\sigma \ov{\zeta}}, \\ 
2 \, y^1_{\zeta} &= U^{\rho}_{\rho \zeta}, & \qquad \qquad
2 \, y^1_{\ov{\zeta}}   &=   U^{\rho}_{\rho \ov{\zeta}}, & \qquad \qquad
-y^2_{\sigma} + y^4_{\rho}  & = U^{\zeta}_{\sigma \rho}, \\
-y^1_{\sigma} + y^4_{\zeta} & = U^{\zeta}_{\sigma \zeta}, & \qquad \qquad
y^4_{\ov{\zeta}} & = U^{\zeta}_{\sigma \ov{\zeta}}, & \qquad \qquad
-y^1_{\rho} + y^2_{\zeta} & = U^{\zeta}_{\rho \zeta}, \\
y^2_{\ov{\zeta}}  & = U^{\zeta}_{\rho \overline{\zeta}}, & \qquad \qquad
 y^1_{\ov{\zeta}} & =  U^{\zeta}_{\zeta \overline{\zeta}}.
\end{alignat*} 
Eliminating $y^1_{\ov{\zeta}}$ among the previous equations leads to:
\begin{equation*}
U^{\zeta}_{\zeta \ov{\zeta}} = \frac{1}{2} \, U^{\rho}_{\rho \ov{\zeta}} = \frac{1}{3} \, U^{\sigma}_{\sigma \ov{\zeta}} 
,\end{equation*}
that is:
\begin{equation*}
{\frac {i\bb}{{\A}^{2}}} = \frac{1}{2} \, \left( {\frac {\cc}{{\A}^{3}}}-{\frac {i\bb}{{\A}^{2}}}  \right)
=- \frac{1}{3} \,  {\frac {\cc}{{\A}^{3}}}
,\end{equation*}
from which we easily deduce that 
\begin{equation*}
\bb = \cc = 0.
\end{equation*}

\subsection{Normalizations of $\dd$ and $\ee$}
We have thus reduced the group $G_2$ to a new group $G_3$, 
whose elements are of the form 
\begin{equation*}
g := \begin{pmatrix}
{\A}^3 & 0 & 0 & 0  \\
0  & \A^2 & 0 & 0 \\
\dd & 0 & \A & 0 \\
\ee & 0 & 0 & \A
\end{pmatrix}.
\end{equation*}

The Maurer Cartan forms of $G_3$ are:
\begin{equation*}
\begin{aligned}
\gamma^1  &: = {\frac {{ d\A}}{\A}}, \\
\gamma^2 & :=-{\frac {\dd { d\A}}{{\A}^{4}}} + {\frac {{
 d\dd}}{{\A}^{3}}}, \\
\gamma^3 & :=-\frac {  \ee d \A}{\A^4} +
\frac { d\ee}{\A^{3}}.
\end{aligned}
\end{equation*}

The third loop of Cartan's method is straightforward.
We get the following structure equations:
\begin{equation*}
\begin{aligned}
d \sigma &= 3 \left. 
\gamma^1 \wedge \sigma \right. 
+ \frac{\dd + \ee}{\A^4} \left. \sigma \wedge \rho \right. + \left. \rho \wedge \zeta \right. 
+ \left.  \rho \wedge \ov{\zeta} \right., \\
d \rho &  =  2 \left. \gamma^1 \wedge \rho \right. + i \, \frac{\ee}{\A^3} \, \left. \sigma \wedge \zeta \right. 
- i \, \frac{\dd}{\A^3} \left. \sigma \wedge \ov{\zeta} \right. + i \, \left. \zeta \wedge \ov{\zeta} \right., \\
d \zeta &= \left. \gamma^1 \wedge \zeta \right. 
+ \left. \gamma^2 \wedge \sigma \right. + \frac{\dd \left( \dd + \ee \right)}{\A^6} 
\left. \sigma \wedge \rho \right.
+
\frac{\dd}{\A^3} \left. \rho \wedge \zeta \right. + \frac{\dd}{\A^3} \left. \rho \wedge \ov{\zeta} \right., \\  
d \ov{\zeta} &= \left. \gamma^1 \wedge \ov{\zeta} \right. + \left. \gamma^3 \wedge \sigma \right. 
+ \frac{\ee \left( \dd + \ee \right)}{\A^6} \left. \sigma \wedge \rho \right.
+
\frac{\ee}{\A^3} \left. \rho \wedge \zeta \right. + \frac{\ee}{\A^3} \left. \rho \wedge \ov{\zeta} \right. 
,\end{aligned}
\end{equation*}
from which we deduce that we can perform the normalizations:
\begin{equation*}
\ee = \dd = 0.
\end{equation*}

With the $1$-dimensional group $G_4$ whose elements $g$ are of the form:
\begin{equation*}
g := \begin{pmatrix}
{\A}^3 & 0 & 0 & 0  \\
0  & \A^2 & 0 & 0 \\
0 & 0 & \A & 0 \\
0 & 0 & 0 & \A
\end{pmatrix},
\end{equation*}
and whose Maurer-Cartan form is given by 
\begin{equation*}
\alpha:= \frac{d \A}{\A},
\end{equation*}
we get the following structure equations:
\begin{equation*}
\begin{aligned}
d \sigma &= 3 \left. \alpha \wedge \sigma \right. + \left. \rho \wedge \zeta \right. 
+ \left.  \rho \wedge \ov{\zeta} \right., \\
d \rho &  =  2 \left. \alpha \wedge \rho \right. + i \, \left. \zeta \wedge \ov{\zeta} \right., \\
d \zeta &= \left. \alpha \wedge \zeta \right., \\
d \ov{\zeta} &= \left. \alpha \wedge \ov{\zeta} \right. 
.\end{aligned}
\end{equation*}
No more normalizations are allowed  at this stage. 
We thus just perform a prolongation by adjoining the form $\alpha$ to the structure equations, whose
exterior derivative is given by:
\begin{equation*}
d \alpha = 0.
\end{equation*}
This completes the proof of Theorem~\ref{thm:Bc}.

\section{Class  ${\sf III_{2}}$}

This section is devoted to the determination of the Lie algebra of CR-automorphisms
of the model manifold of class ${\sf III_{2}}$ which is defined by the equations:
\begin{equation*}
{\sf N}: \qquad \qquad
\begin{aligned}
w_1 & = \ov{w_1} + 2 \, i \, z \ov{z}, \\
w_2 & = \ov{w_2} + 2 \, i \, z \ov{z} \left( z + \ov{z} \right), \\
w_3 & = \ov{w_3} + 2 i \, z \zb ( z^2 + \frac{3}{2} \, z \zb + \zb^2 ).
\end{aligned}
\end{equation*}
It is the model manifold 
for CR-manifolds belonging to class ${\sf III_2}$, 
in the sense that any such manifold might be viewed
as a deformation of ${\sf N}$ by the way of a Cartan connection ( \cite{pocchiola3}).
The main result of this section is the following:

\begin{theorem}
\label{thm:N}
The model of the class ${\sf III_2}${\rm :}
\begin{equation*}
{\sf N}: \qquad \qquad
\begin{aligned}
w_1 & = \ov{w_1} + 2 \, i \, z \ov{z}, \\
w_2 & = \ov{w_2} + 2 \, i \, z \ov{z} \left( z + \ov{z} \right), \\
w_3 & = \ov{w_3} + 2 i \, z \zb ( z^2 + \frac{3}{2} \, z \zb + \zb^2 ),
\end{aligned}
\end{equation*}
has a ${\bf 6}$-dimensional Lie algebra of CR-automorphisms. 
A basis for the Maurer-Cartan forms of ${\sf aut_{CR}}({\sf N})$ is provided
by the $6$ differential $1$-forms  $\tau$, $\sigma$, $\rho$, $\zeta$,  $\ov{\zeta}$, $\alpha$,
which satisfy the structure equations:
\begin{equation*}
\begin{aligned}
d \tau &= 4 \left. \alpha \wedge \tau \right.
+ \left.\sigma \wedge \zeta \right.
+ \left. \sigma \wedge \ov{\zeta} \right., \\
d \sigma &= 3 \left. \alpha \wedge \sigma \right. + \left. \rho \wedge \zeta \right. 
+ \left.  \rho \wedge \ov{\zeta} \right., \\
d \rho &  =  2 \left. \alpha \wedge \rho \right. + i \, \left. \zeta \wedge \ov{\zeta} \right., \\
d \zeta &= \left. \alpha \wedge \zeta \right., \\
d \ov{\zeta} &= \left. \alpha \wedge \ov{\zeta} \right.,\\
d \alpha &= 0.
\end{aligned}
\end{equation*}

\end{theorem}

\subsection{Initial $G$-structure}

The vector fields :
\[
\mathcal{L}
:=
\frac{\partial}{\partial z}
+
i\overline{z}\,\frac{\partial}{\partial u_1}
+
i\big(2z\overline{z}+\overline{z}^2\big)\,
\frac{\partial}{\partial u_2}
+
i\big(3z^2\overline{z}
+
3z\overline{z}^2
+
\overline{z}^3
\big)\,
\frac{\partial}{\partial u_3},
\]
with its conjugate:
\[
\overline{\mathcal{L}}
:=
\frac{\partial}{\partial\overline{z}}
-
iz\,\frac{\partial}{\partial u_1}
-
i\big(2z\overline{z}+z^2\big)\,
\frac{\partial}{\partial u_2}
-
i\big(
3z\overline{z}^2
+
3z^2\overline{z}
+
z^3
\big)\,
\frac{\partial}{\partial u_3},
\]
constitute a basis of $T^{1,0}_p {\sf N}$ and of $T^{0,1}_p {\sf N}$ at each point $p$ of ${\sf N}$.
Moreover the vector fields $\T$, $\s$ and $\mathcal{R}$ defined by:
\begin{equation*}
\T: = i \, \big[\mathcal{L}, \Lb] 
,\end{equation*}

\begin{equation*}
\s: = \big[\LL, \T \big], 
\end{equation*}
and

\begin{equation*}
\mathcal{R} := \big[\LL, \s \big],
\end{equation*}
complete a frame on ${\sf N}$:
\begin{equation*}
\big\{\mathcal{R}, \mathcal{S},\,
\mathcal{T},\,\mathcal{L}, \,
\overline{\mathcal{L}}
\big\}
.
\end{equation*}
The expressions of $\T$, $\s$ and $\mathcal{R}$ are:

\begin{equation*}
\begin{aligned}
\mathcal{T}
:=
&
\,\,2\,\frac{\partial}{\partial u_1}
+
\big(4z+4\overline{z}\big)
\frac{\partial}{\partial u_2}
+
\big(6z^2+12z\overline{z}+6\overline{z}^2\big)\,
\frac{\partial}{\partial u_3}, \\
\s
:=
&
\,\,4\,\frac{\partial}{\partial u_2}
+
\big(12z+12\overline{z}\big)\,
\frac{\partial}{\partial u_3},
\\
\mathcal{R}
:=
&
\,\,12\,\frac{\partial}{\partial u_3}.
\end{aligned}
\end{equation*}
The dual coframe $ \big\{ \tau_0,\,\sigma_0,\,
\rho_0,\,\overline{\zeta_0},\,
\zeta_0
\big\}$
is thus given by:

\begin{equation*}
\begin{aligned}
\tau_0 &=  - \frac{i}{12} \, \zb^3 \, dz 
+ \frac{i}{12} \, z^3 \, d \zb + \left( \frac{1}{4} 
\, z^2 + \frac{1}{2} \, z \zb + \frac{1}{4} \, \zb^2 \right) du_1
- \left( \frac{1}{4} \, z + \frac{1}{4} \, \zb \right) du_2 + \frac{1}{12} \, du_3, \\
\sigma_0 &= \frac{i}{4} \, \zb^2 \, dz - 
\frac{i}{4} \, z^2 \, d \zb -  \left( \frac{1}{2} 
\, z + \frac{1}{2} \, \zb \right) du_1 + \frac{1}{4} \, du_2, \\
\rho_0 &=  - \frac{i}{2} \, \zb \,  dz + \frac{i}{2} \, z \, d\zb + \frac{1}{2} \, du_1 ,\\
\zeta_0 &= dz, \\
\ov{\zeta_0} &= d\zb.
\end{aligned}
\end{equation*}

We deduce the structure equations enjoyed by the base coframe
$
\big\{ \tau_0,\,\sigma_0,\,
\rho_0,\,\overline{\zeta_0},\,
\zeta_0
\big\}
$:
\begin{equation}
\label{eq:N}
\begin{aligned}
d \tau_{0} & =   \sigma_0 \wedge \zeta_0 + \sigma_0 \wedge \ov{\zeta_0}, \\
d \sigma_0 & = \rho_0 \wedge \zeta_0 + \rho_0 \wedge \ov{\zeta_0}, \\
d \rho_{0} & =  i \, \zeta_0 \wedge \ov{\zeta_0},\\
d \zeta_0 & = 0, \\
d \ov{ \zeta_{0}} & = 0.
\end{aligned}
\end{equation}

As the torsion coefficients of these structure equations are constants, we have the following result:

\begin{lemma}
The model of the class ${\sf III_2}$
is locally isomorphic to a Lie group whose Maurer-Cartan forms satisfy the structure equations 
$(\ref{eq:N})$.
\end{lemma} 

The matrix Lie group which encodes suitably the equivalence problem for the model of class ${\sf III_2}$ (see \cite{pocchiola3}) is the 
$18$-dimensional Lie group $G_1$ whose elements $g$ are of the form:
\begin{equation*}
g := \begin{pmatrix}
{\A^3} \ab & 0 & 0 & 0 & 0 \\
\ff  & \A^2 \ab & 0 & 0 & 0 \\
\G & \cc & \A \ab & 0 & 0 \\
\hh & \dd & \bb & \A & 0 \\
\kk & \ee & \bbb & 0 & \ab
\end{pmatrix}.
\end{equation*}

With the notations:

\begin{alignat*}{1}
\omega_{0} :=
\begin{pmatrix}
\tau_0 \\ \sigma_0 \\ \rho_0 \\ \zeta_{0} \\  \ov{\zeta}_0
\end{pmatrix}, & \qquad
\omega:=
\begin{pmatrix}
\tau \\ \sigma \\ \rho \\ \zeta \\ \ov{\zeta}
\end{pmatrix},
\end{alignat*}
we introduce the $G_1$-structure $P^1$ on ${\sf N}$ constituted 
by the coframes $\omega$ which satisfy the relation:

\begin{equation*}
\omega:= g \cdot \omega_0.
\end{equation*}
As in the case of Beloshapka's cubic,
the proof of theorem (\ref{thm:N}) relies on successive reductions of $P^1$ through Cartan's equivalence method.

\subsection{Normalization of {\A}}

The Maurer-Cartan forms of $G_1$ are given by:
\begin{equation*}
\begin{aligned}
\alpha^1  &: = {\frac {{ d\A}}{\A}}, \\
\alpha^2 &:= -{\frac {\bb{ d\A}}{{\A}^{2}{ \ab}}}+{\frac {{ d\bb}}{\A{ \ab}}}, \\
\alpha^3  &:=-{\frac {\cc{ d\A}}{{ \ab}\,{\A}^{3}}}-{\frac {\cc{ d\ab
}}{{{ \ab}}^{2}{\A}^{2}}} + { \frac {{ d\cc}}{{\A}^{2}{ \ab}}}, \\
\alpha^4 & :=-{\frac { \left( \dd\A{ \ab}-\bb\cc \right) { d\A}}{{\A}^{4}{
{ \ab}}^{2}}}-{\frac {\cc{ d\bb}}{{\A}^{3}{{ \ab}}^{2}}}+{\frac {{
 d\dd}}{{\A}^{2}{ \ab}}}, \\
\alpha^5 & :=-{\frac { \left( \ee\A{ \ab}-{ \bbb}\,\cc \right) { d\ab}
}{{\A}^{3}{{ \ab}}^{3}}}-{\frac {\cc{ d\bbb}}{{\A}^{3}{{ \ab}}^{2}}}+
{\frac {{ d\ee}}{{\A}^{2}{ \ab}}}, \\
\alpha^6 & :=-2\,{\frac {\ff{ d\A}}{{ \ab}\,{\A}^{4}}}-{\frac {\ff{ 
d\ab}}{{\A}^{3}{{ \ab}}^{2}}}+{\frac {{ d\ff}}{{ \ab}\,{\A}^{3}}}, \\
\alpha^7 & :=-{\frac { \left( \G{\A}^{2}{ \ab}-\cc\ff \right) { d\A}}{{{
 \ab}}^{2}{\A}^{6}}}-{\frac { \left( \G{\A}^{2}{ \ab}-\cc\ff \right) {
 d\ab}}{{{ \ab}}^{3}{\A}^{5}}}-{\frac {\ff{ d\cc}}{{\A}^{5}{{ \ab}}
^{2}}}+{\frac {{ d\G}}{{ \ab}\,{\A}^{3}}}, \\
\alpha^8 & :=-{\frac { \left( \hh{\A}^{3}{{ \ab}}^{2}- \dd \ff\A{ \ab}-\bb \G{\A}
^{2}{ \ab}+\bb\cc\ff \right) { d\A}}{{\A}^{7}{{ \ab}}^{3}}}-{\frac {
 \left( \G{\A}^{2}{ \ab}-\cc\ff \right) { d\bb}}{{\A}^{6}{{ \ab}}^{3}}}-
{\frac {\ff{ d\dd}}{{\A}^{5}{{ \ab}}^{2}}}+{\frac {{ d\hh}}{{ \ab}
\,{\A}^{3}}},\\
\alpha^9 & :=-{\frac { \left( \kk{\A}^{3}{{ \ab}}^{2}-\ee\ff\A{ \ab}-{ 
\bbb}\,\G{\A}^{2}{ \ab}+{ \bbb}\,\cc\ff \right) { d\ab}}{{\A}^{6}{{ \ab}
}^{4}}}-{\frac { \left( \G{\A}^{2}{ \ab}-\cc\ff \right) { d\bbb}}{{\A}^{6}
{{ \ab}}^{3}}}-{\frac {\ff{ d\ee}}{{\A}^{5}{{ \ab}}^{2}}}+{\frac {{
 d\kk}}{{ \ab}\,{\A}^{3}}},
\end{aligned}
\end{equation*}
together with their conjugates.

The first structure equation is given by:
\begin{multline*}
d \tau = 3 \left.\alpha^1 \wedge \tau \right. + \left.\ov{\alpha^1} \wedge \tau \right. \\
+ T^{\tau}_{\tau \sigma} \left. \tau \wedge \sigma \right.
+ T^{\tau}_{\tau \rho} \left. \tau \wedge \rho \right.
+ T^{\tau}_{\tau \zeta} \left. \tau \wedge \zeta \right. \\
+ T^{\tau}_{\tau \ov{\zeta}} \left. \tau \wedge \ov{\zeta} \right.
+ T^{\tau}_{\sigma \rho} \left.\sigma \wedge \rho \right.
+ \left.\sigma \wedge \zeta \right.
-  \frac{\A}{\ab} \, \left. \sigma \wedge \ov{\zeta} \right. 
,\end{multline*}
from which we immediately 
deduce that $\frac{\A}{\ab}$ is an essential torsion coefficient which shall be normalized to $1$ by setting:
\begin{equation*}
\A = \ab.
\end{equation*}

We thus have reduced the $G_{1}$ equivalence 
problem to a $G_2$ equivalence problem, where $G_2$ is the $10$ dimensional real matrix 
Lie group whose elements are of the form
\begin{equation*}
g = \begin{pmatrix} 
{\A^4}  & 0 & 0 & 0 & 0 \\
\ff  & \A^3 & 0 & 0 & 0 \\
\G & \cc & \A^2 & 0 & 0 \\
\hh & \dd & \bb & \A & 0 \\
\kk & \ee & \bbb & 0 & \A
\end{pmatrix}
,\end{equation*}

\subsection{Normalizations of $\ff$, $\bb$ and $\cc$}

The Maurer-Cartan forms of $G_2$ are given by:
\begin{equation*}
\begin{aligned}
\beta^1  &: = {\frac {{ d\A}}{\A}}, \\
\beta^2 &:= -{\frac {\bb d\A }{{\A}^{3}}}+{\frac {{ d\bb}}{\A^2}}, \\
\beta^3  &:=- 2 \, {\frac {\cc  d\A }{{\A}^{4}}} + { \frac {{ d\cc}}{{\A}^{3}}}, \\
\beta^4 &: =-{\frac { \left( \dd\A^2 -\bb \cc \right) { d\A}}{{\A}^{6}}}-{\frac {\cc{ d\bb}}{{\A}^{5}}}+{\frac {{
 d\dd}}{{\A}^{3}}}, \\
\beta^5 & :=-\frac { \left( \ee\A^2-{ \bbb}\,\cc \right)  d\A
}{\A^6}-\frac {\cc{ d\bbb}}{\A^{5}}+
\frac { d\ee}{\A^{3}}, \\
\beta^6 & := -3 \,{\frac {\ff{ d\A}}{{\A}^{5}}} + {\frac {{ d\ff}}{{\A}^{4}}}, \\
\beta^7 & :=- 2 \, {\frac { \left( \G{\A}^{3}-\cc \ff \right) 
{ d\A}}{\A^{8}}} -{\frac {\ff{ d\cc}}{{\A}^{7}}}+{\frac {{ d\G}}{{\A}^{4}}}, \\
\beta^8 & :=-{\frac { \left( \hh{\A}^{5}- \dd \ff\A^2-\bb \G{\A}
^{3}+\bb\cc\ff \right) { d\A}}{{\A}^{10}}}-{\frac {
 \left( \G{\A}^{3}-\cc\ff \right) { d\bb}}{{\A}^{9}}}-
{\frac {\ff{ d\dd}}{\A^{7}}}+{\frac {{ d\hh}}{
{\A}^{4}}},\\
\beta^9 & :=-{\frac { \left( \kk{\A}^{5}-\ee\ff\A^2-{ 
\bbb}\,\G{\A}^{3}+{ \bbb}\,\cc\ff \right) { d\A}}{{\A}^{10}}}-{\frac { \left( \G{\A}^{3}-\cc\ff \right) { d\bbb}}{{\A}^{9}}}
-
{\frac {\ff{ d\ee}}{{\A}^{5}{{ \ab}}^{2}}}+{\frac {{
 d\kk}}{{\A}^{4}}},
\end{aligned}
\end{equation*}
together with $\ov{\beta^i}, \,\,\, i=2 \dots 9$.

Using formula (\ref{eq:str}), we get the structure 
equations for the lifted coframe $(\tau, \sigma, \rho, \zeta, \ov{\zeta})$ from those of the base coframe
$(\tau_0, \sigma_0, \rho_0, \hat{\zeta}_0, \ov{\hat{\zeta}_0})$
by a matrix multiplication and a linear change of coordinates, as in the first  step:

\begin{multline*}
d \tau 
=
4 \, \beta^1 \wedge \tau \\
+
U^{\tau}_{\tau \sigma} \, \tau \wedge \sigma
+
U^{\tau}_{\tau \rho} \, \tau \wedge \rho
+
U^{\tau}_{\tau \zeta} \, \tau \wedge \zeta 
+
U^{\tau}_{\tau \ov{\zeta}} \, \tau \wedge \ov{\zeta} \\
+
U^{\tau}_{\sigma \rho} \, \sigma \wedge \rho
+
\sigma \wedge \zeta
+
\sigma \wedge \ov{\zeta},
\end{multline*}

\begin{multline*}
d \sigma = 3 \, \beta^1 \wedge \sigma + \beta^6 \wedge \tau \\
+
U^{\sigma}_{\tau \sigma} \left. \tau \wedge \sigma \right.
+
U^{\sigma}_{\tau \rho} \left. \tau \wedge \rho \right.
+
U^{\sigma}_{\tau \zeta} \left. \tau \wedge \zeta \right. \\
+
U^{\sigma}_{\tau \ov{\zeta}} \left. \tau \wedge \ov{\zeta} \right.
+
U^{\sigma}_{\sigma \rho} \left. \sigma \wedge \rho \right.
+
U^{\sigma}_{\sigma \zeta} \left. \sigma \wedge \zeta \right. \\
+
U^{\sigma}_{\sigma \ov{\zeta}} \left. \sigma \wedge \ov{\zeta} \right.
+
\rho \wedge \zeta
+
\rho \wedge \ov{\zeta},
\end{multline*}

\begin{multline*}
d \rho
=
2 \beta^{1} \wedge \rho + \beta^3 \wedge \sigma + \beta^7 \wedge \tau \\
+
U^{\rho}_{\tau \sigma} \, \tau \wedge \sigma
+
U^{\rho}_{\tau \rho} \, \tau \wedge \rho
+
U^{\rho}_{\tau \zeta} \, \tau \wedge \zeta 
+
U^{\rho}_{\tau \overline{\zeta}} \, \rho \wedge \overline{\zeta} 
+
U^{\rho}_{\sigma \rho} \, \sigma \wedge \rho \\
+
U^{\rho}_{\sigma \zeta} \, \sigma \wedge \zeta
+
U^{\rho}_{\sigma \ov{\zeta}} \, \sigma \wedge \ov{\zeta} 
+
U^{\rho}_{\rho \zeta} \, \rho \wedge \zeta  
+
U^{\rho}_{\rho \ov{\zeta}} \, \rho \wedge \ov{\zeta} 
+
i \, \zeta \wedge \ov{\zeta}
,\end{multline*}

\begin{multline*}
d \zeta
=
{\beta}^{1} \wedge \zeta + {\beta}^{2} \wedge \rho + {\beta}^{4} \wedge \sigma + \beta^8 \wedge \tau \\
+
U^{\zeta}_{\tau \sigma} \, \tau \wedge \sigma
+
U^{\zeta}_{\tau \rho} \, \tau \wedge \rho
+
U^{\zeta}_{\tau \zeta} \, \tau \wedge \zeta 
+
U^{\zeta}_{\tau \ov{\zeta}} \, \tau \wedge \ov{\zeta} \\
+
U^{\zeta}_{\sigma \rho} \, \sigma \wedge \rho
+
U^{\zeta}_{\sigma \zeta} \, \sigma \wedge \zeta 
+
U^{\zeta}_{\sigma \ov{\zeta}} \, \sigma \wedge \ov{\zeta} 
+
U^{\zeta}_{\rho \zeta} \, \rho \wedge \zeta \\
+
U^{\zeta}_{\rho \overline{\zeta}} \, \rho \wedge \overline{\zeta}
+
U^{\zeta}_{\zeta \overline{\zeta}} \, \zeta \wedge \overline{\zeta}
.\end{multline*}

We now proceed with the absorption phase. We introduce the modified Maurer-Cartan forms:
\begin{equation*}
\widetilde{\beta}^i= \beta^i  -y_{\tau}^i \, \tau - y_{\sigma} \, \sigma - y_{\rho}^i \, \rho  - y_{\zeta}^i \, \zeta \, 
- y_{\overline{\zeta}}^i \, \overline{\zeta}.
\end{equation*}
The structure equations rewrite:
\begin{multline*}
d \tau 
=
4 \, \widetilde{\beta}^1 \wedge \tau \\
+
\left( U^{\tau}_{\tau \sigma} - 4 \, y^1_{\sigma} \right) \, \left. \tau \wedge \sigma \right.
+
\left( U^{\tau}_{\tau \rho} - 4 \, y^1_{\rho} \right) \, \left. \tau \wedge \rho \right. \\
+
\left(U^{\tau}_{\tau \zeta} - 4 \, y^1_{\zeta} \right) \, \left. \tau \wedge \zeta \right. 
+
\left(U^{\tau}_{\tau \ov{\zeta}} - 4 \,  y^1_{\ov{\zeta}} \right) \, \left. \tau \wedge \ov{\zeta} \right. \\
+
U^{\tau}_{\sigma \rho} \, \left. \sigma \wedge \rho \right.
+
\left.\sigma \wedge \zeta \right.
+
\left. \sigma \wedge \ov{\zeta} \right.,
\end{multline*}

\begin{multline*}
d \sigma = 3 \, \left. \widetilde{\beta}^1 \wedge \sigma \right. + \left. \widetilde{\beta}^6 \wedge \tau \right. \\
+
\left( U^{\sigma}_{\tau \sigma} + 3 \, y^1_{\tau} - y^6_{\sigma} \right) \, \left. \tau \wedge \sigma \right.
+
\left( U^{\sigma}_{\tau \rho} -  y^6_{\rho} \right) \, \left. \tau \wedge \rho \right. \\
+
\left( U^{\sigma}_{\tau \zeta} - y^6_{\zeta} \right) \, \left. \tau \wedge \zeta \right. 
+
\left( U^{\sigma}_{\tau \ov{\zeta}} - y^6_{\ov{\zeta}} \right) \, \left. \tau \wedge \ov{\zeta} \right. \\
+
\left( U^{\sigma}_{\sigma \rho} - 3 \, y^1_{\rho} \right) \, \left. \sigma \wedge \rho \right. 
+
\left( U^{\sigma}_{\sigma \zeta} -3 \,y^1_{\zeta} \right) \, \left. \sigma \wedge \zeta \right. \\
+
\left( U^{\sigma}_{\sigma \ov{\zeta}}  -3 \,y^1_{\ov{\zeta}} \right) \, \left. \sigma \wedge \ov{\zeta} \right.
+
\rho \wedge \zeta
+
\rho \wedge \ov{\zeta},
\end{multline*}

\begin{multline*}
d \rho
=
2 \widetilde{\beta}^{1} \wedge \rho + \widetilde{\beta}^3 \wedge \sigma + \widetilde{\beta}^7 \wedge \tau \\
+
\left( U^{\rho}_{\tau \sigma} + y^3_{\tau} - y^7_{\sigma} \right) \, \left. \tau \wedge \sigma \right.
+
\left( U^{\rho}_{\tau \rho} + 2 \, y^1_{\tau} - y^7_{\rho} \right) \, \left. \tau \wedge \rho \right. \\
+
\left( U^{\rho}_{\tau \zeta} - y^7_{\zeta} \right) \, \left. \tau \wedge \zeta \right. 
+
\left( U^{\rho}_{\tau \overline{\zeta}} - y^7_{\ov{\zeta}} \right) \, \left. \rho \wedge \overline{\zeta} \right. \\
+
\left( U^{\rho}_{\sigma \rho} + 2 \, y^1_{\sigma} - y^3_{\rho} \right) \, \left. \sigma \wedge \rho \right.
+
\left( U^{\rho}_{\sigma \zeta} - y^3_{\zeta} \right) \, \left. \sigma \wedge \zeta \right. \\
+
\left( U^{\rho}_{\sigma \ov{\zeta}} - y^3_{\ov{\zeta}} \right) \, \left. \sigma \wedge \ov{\zeta} \right.
+
\left( U^{\rho}_{\rho \zeta} - 2 \, y^1_{\zeta} \right) \, \left. \rho \wedge \zeta \right.  \\ 
+
\left( U^{\rho}_{\rho \ov{\zeta}} - 2 \, y^1_{\ov{\zeta}} \right) \, \left. \rho \wedge \ov{\zeta} \right.
+
i \, \left. \zeta \wedge \ov{\zeta} \right.
,\end{multline*}

\begin{multline*}
d \zeta
=
\widetilde{\beta^{1}} \wedge \zeta +
 \widetilde{\beta^{2}} \wedge \rho + \widetilde{\beta^{4}} \wedge \sigma + \widetilde{\beta^8} \wedge \tau \\
+
\left( U^{\zeta}_{\tau \sigma} + y^4_{\tau}- y^8_{\sigma} \right) \, \left. \tau \wedge \sigma \right.
+
\left( U^{\zeta}_{\tau \rho} + y^2_{\tau} - y^8_{\rho} \right) \, \left. \tau \wedge \rho \right. \\
+
\left( U^{\zeta}_{\tau \zeta} + y^1_{\tau} - y^8_{\zeta} \right) \, \left. \tau \wedge \zeta \right.
+
\left( U^{\zeta}_{\tau \ov{\zeta}} - y^8_{\ov{\zeta}} \right) \, \left. \tau \wedge \ov{\zeta} \right. \\
+
\left( U^{\zeta}_{\sigma \rho} + y^2_{\sigma}- y^4_{\rho} \right) \, \left. \sigma \wedge \rho \right.
+
\left( U^{\zeta}_{\sigma \zeta} + y^1_{\sigma} - y^4_{\zeta} \right) \, \left. \sigma \wedge \zeta \right. \\ 
+
\left( U^{\zeta}_{\sigma \ov{\zeta}} - y^4_{\ov{\zeta}} \right) \, \left. \sigma \wedge \ov{\zeta} \right. 
+
\left( U^{\zeta}_{\rho \zeta} + y^1_{\rho} - y^2_{\zeta} \right) \, \left. \rho \wedge \zeta \right. \\
+
\left( U^{\zeta}_{\rho \overline{\zeta}} - y^2_{\ov{\zeta}} \right) \, \left. \rho \wedge \overline{\zeta} \right.
+
\left( U^{\zeta}_{\zeta \overline{\zeta}} - y^1_{\ov{\zeta}} \right) \, \left. \zeta \wedge \overline{\zeta} \right.
.\end{multline*}

We get the following absorbtion equations:

\begin{alignat*}{3}
4 \, y^1_{\sigma} &=U^{\tau}_{\tau \sigma},
& \qquad \qquad  4 \, y^1_{\rho}  &=  U^{\tau}_{\tau \rho} ,
& \qquad \qquad  4 \, y^1_{\zeta} &=  U^{\tau}_{\tau \zeta} , \\
4 \,  y^1_{\ov{\zeta}} &= U^{\tau}_{\tau \ov{\zeta}}, 
& \qquad \qquad - 3 \, y^1_{\tau} + y^6_{\sigma}  & =U^{\sigma}_{\tau \sigma},  & 
\qquad \qquad y^6_{\rho}   & =  U^{\sigma}_{\tau \rho} , \\
y^6_{\zeta}  &= U^{\sigma}_{\tau \zeta},  & \qquad \qquad  
y^6_{\ov{\zeta}}  & =   U^{\sigma}_{\tau \ov{\zeta}}, & \qquad \qquad 
3 \, y^1_{\rho} & =  U^{\sigma}_{\sigma \rho}, \\  
3 \,y^1_{\zeta} & =  U^{\sigma}_{\sigma \zeta}, & \qquad \qquad
3 \,y^1_{\ov{\zeta}} & =  U^{\sigma}_{\sigma \ov{\zeta}}, & \qquad \qquad
-y^3_{\tau} + y^7_{\sigma}  & = U^{\rho}_{\tau \sigma},  \\ 
-2 \, y^1_{\tau} + y^7_{\rho} &=  U^{\rho}_{\tau \rho}, & \qquad \qquad 
y^7_{\zeta}  & = U^{\rho}_{\tau \zeta}, &\qquad \qquad   
y^7_{\ov{\zeta}}  &=  U^{\rho}_{\tau \overline{\zeta}}, \\
-2 \, y^1_{\sigma} + y^3_{\rho}  &= U^{\rho}_{\sigma \rho}, &\qquad \qquad 
y^3_{\zeta}  &=  U^{\rho}_{\sigma \zeta}  , &\qquad \qquad   
 y^3_{\ov{\zeta}}   &= U^{\rho}_{\sigma \ov{\zeta}}, \\ 
2 \, y^1_{\zeta} &= U^{\rho}_{\rho \zeta}, & \qquad \qquad
2 \, y^1_{\ov{\zeta}}   &=   U^{\rho}_{\rho \ov{\zeta}}, & \qquad \qquad
-y^4_{\tau} + y^8_{\sigma} & = U^{\zeta}_{\tau \sigma},\\
-y^2_{\tau} + y^8_{\rho} & =  U^{\zeta}_{\tau \rho},& \qquad \qquad 
- y^1_{\tau} + y^8_{\zeta} & = U^{\zeta}_{\tau \zeta} , & \qquad \qquad
y^8_{\ov{\zeta}} & =U^{\zeta}_{\tau \ov{\zeta}} , \\
-y^2_{\sigma} + y^4_{\rho}  & = U^{\zeta}_{\sigma \rho}, & \qquad \qquad 
-y^1_{\sigma} + y^4_{\zeta} & = U^{\zeta}_{\sigma \zeta}, & \qquad \qquad
y^4_{\ov{\zeta}} & = U^{\zeta}_{\sigma \ov{\zeta}}, \\
-y^1_{\rho} + y^2_{\zeta} & = U^{\zeta}_{\rho \zeta}, & \qquad \qquad 
y^2_{\ov{\zeta}}  & = U^{\zeta}_{\rho \overline{\zeta}}, & \qquad \qquad
 y^1_{\ov{\zeta}} & =  U^{\zeta}_{\zeta \overline{\zeta}}.
\end{alignat*} 
Eliminating $y^1_{\ov{\zeta}}$ among the previous equations leads to:
\begin{equation*}
U^{\zeta}_{\zeta \ov{\zeta}} = \frac{1}{2} \, U^{\rho}_{\rho \ov{\zeta}} = 
\frac{1}{3} \, U^{\sigma}_{\sigma \ov{\zeta}} = \frac{1}{4} \, U^{\tau}_{\tau \ov{\zeta}}
,\end{equation*}
that is:
\begin{equation*}
{\frac {i\bb}{{\A}^{2}}} = \frac{1}{2} \, \left( {\frac {\cc}{{\A}^{3}}}-{\frac {i\bb}{{\A}^{2}}}  \right)
=- \frac{1}{3} \,\left(  {\frac {\cc}{{\A}^{3}}}+{\frac {\ff}{{\A}^{4}}} \right) =-\frac{1}{4} \,{\frac {\ff}{{\A}^{4}}},
\end{equation*}
from which we easily deduce that 
\begin{equation*}
\bb = \cc = \ff = 0.
\end{equation*}

We have thus reduced the group $G_2$ to a new group $G_3$, 
whose elements are of the form 
\begin{equation*}
g := \begin{pmatrix}
{\A^4}  & 0 & 0 & 0 & 0 \\
0  & \A^3 & 0 & 0 & 0 \\
\G & 0 & \A^2 & 0 & 0 \\
\hh & \dd & 0 & \A & 0 \\
\kk & \ee & 0 & 0 & \A
\end{pmatrix}.
\end{equation*}

\subsection{Normalization of $\G$, $\dd$ and $\ee$}
The Maurer Cartan forms of $G_3$ are:
\begin{equation*}
\begin{aligned}
\gamma^1  &: = {\frac {{ d\A}}{\A}}, \\
\gamma^2 & :=-{\frac {\dd { d\A}}{{\A}^{4}}} + {\frac {{
 d\dd}}{{\A}^{3}}}, \\
\gamma^3 &: =-\frac {  \ee d \A}{\A^4} +
\frac { d\ee}{\A^{3}}, \\
\gamma^4 &: =- 2 \, \frac {  \G d \A}{\A^5} +
\frac { d\G}{\A^{4}}, \\
\gamma^5 &: =- \frac {  \hh d \A}{\A^5} +
\frac { d\hh}{\A^{4}}, \\
\gamma^6 &: =- \frac {  \kk d \A}{\A^5} +
\frac { d\kk}{\A^{4}}.
\end{aligned}
\end{equation*}

We get the following structure equations:
\begin{equation*}
d \tau = 4 \left. \gamma^1 \wedge \tau \right. 
+ V^{\tau}_{\tau \sigma} \, \left. \tau \wedge \sigma \right. + \sigma \wedge \zeta + \sigma \wedge \ov{\zeta},
\end{equation*}

\begin{equation*}
d \sigma = 3 \left. \gamma^1 \wedge \sigma \right. 
+
V^{\sigma}_{\tau \rho} \left. \tau \wedge \rho \right.
+
V^{\sigma}_{\tau \zeta} \left. \tau \wedge \zeta \right.
+
V^{\sigma}_{\tau \ov{\zeta}} \left. \tau \wedge \ov{\zeta} \right.
+
V^{\sigma}_{\sigma \rho} \left. \sigma \wedge \rho \right.
+
\left. \rho \wedge \zeta \right.
+
\left. \rho \wedge \ov{\zeta} \right.
,\end{equation*}

\begin{equation*}
d \rho
=
2 \left. \gamma^{1} \wedge \rho \right. + \left. \gamma^4 \wedge \tau \right.
+
V^{\rho}_{\tau \sigma} \left. \tau \wedge \sigma \right.
+
V^{\rho}_{\tau \zeta}  \left. \tau \wedge \zeta \right.
+
V^{\rho}_{\tau \overline{\zeta}} \left.  \rho \wedge \overline{\zeta} \right. 
+
V^{\rho}_{\sigma \zeta} \left. \sigma \wedge \zeta \right.
+
V^{\rho}_{\sigma \ov{\zeta}} \left. \sigma \wedge \ov{\zeta} \right. 
+
i \, \left. \zeta \wedge \ov{\zeta} \right.
,\end{equation*}

\begin{multline*}
d \zeta
=
\left. \gamma^{1} \wedge \zeta \right.+ \left. {\gamma}^{2} \wedge \sigma \right. 
+ \left. \gamma^5 \wedge \tau \right. \\
+
V^{\zeta}_{\tau \sigma} \, \tau \wedge \sigma
+
V^{\zeta}_{\tau \rho} \, \tau \wedge \rho
+
V^{\zeta}_{\tau \zeta} \, \tau \wedge \zeta 
+
V^{\zeta}_{\tau \ov{\zeta}} \, \tau \wedge \ov{\zeta} \\
+
V^{\zeta}_{\sigma \rho} \, \sigma \wedge \rho
+
V^{\zeta}_{\sigma \zeta} \, \sigma \wedge \zeta 
+
V^{\zeta}_{\sigma \ov{\zeta}} \, \sigma \wedge \ov{\zeta} 
+
V^{\zeta}_{\rho \zeta} \, \rho \wedge \zeta \\
+
V^{\zeta}_{\rho \overline{\zeta}} \, \rho \wedge \overline{\zeta}
,\end{multline*}
and 
\begin{multline*}
d \ov{\zeta}
=
\left. \gamma^{1} \wedge \ov{\zeta} \right.+ \left. \gamma^{3} \wedge \sigma \right. 
+ \left. \gamma^6 \wedge \tau \right. \\
+
V^{\ov{\zeta}}_{\tau \sigma} \, \tau \wedge \sigma
+
V^{\ov{\zeta}}_{\tau \rho} \, \tau \wedge \rho
+
V^{\ov{\zeta}}_{\tau \zeta} \, \tau \wedge \zeta 
+
V^{\ov{\zeta}}_{\tau \ov{\zeta}} \, \tau \wedge \ov{\zeta} \\
+
V^{\ov{\zeta}}_{\sigma \rho} \, \sigma \wedge \rho
+
V^{\ov{\zeta}}_{\sigma \zeta} \, \sigma \wedge \zeta 
+
V^{\ov{\zeta}}_{\sigma \ov{\zeta}} \, \sigma \wedge \ov{\zeta} 
+
V^{\ov{\zeta}}_{\rho \zeta} \, \rho \wedge \zeta \\
+
V^{\ov{\zeta}}_{\rho \overline{\zeta}} \, \rho \wedge \overline{\zeta}
.\end{multline*}

From these equations, we immediately see that $V^{\sigma}_{\tau \zeta}$, 
$V^{\zeta}_{\rho \ov{\zeta}}$ and $V^{\ov{\zeta}}_{\rho \zeta}$ are essential torsion coefficients.
As we have:
\begin{alignat*}{3}
V^{\sigma}_{\tau \zeta} &= - \frac{\G}{\A^4}, \qquad \qquad & V^{\zeta}_{\rho \ov{\zeta}} & = \frac{\dd}{\A^3} , \qquad \qquad & 
V^{\ov{\zeta}}_{\rho \zeta} &=\frac{\ee}{\A^3},
\end{alignat*}
we obtain the new normalizations:
\begin{equation*}
\dd = \ee = \G = 0.
\end{equation*}

The reduced group $G_4$ is of the form:
\begin{equation*}
g := \begin{pmatrix}
{\A^4}  & 0 & 0 & 0 & 0 \\
0  & \A^3 & 0 & 0 & 0 \\
0 & 0 & \A^2 & 0 & 0 \\
\hh & 0 & 0 & \A & 0 \\
\kk & 0 & 0 & 0 & \A
\end{pmatrix}.
\end{equation*}
Its Maurer-Cartan forms are given by:
\begin{equation*}
\begin{aligned}
\delta^1  &: = {\frac {{ d\A}}{\A}}, \\
\delta^2 &: =- \frac {  \hh d \A}{\A^5} +
\frac { d\hh}{\A^{4}}, \\
\delta^3 & :=- \frac {  \kk d \A}{\A^5} +
\frac { d\kk}{\A^{4}}.
\end{aligned}
\end{equation*}
The structure equations are easily computed as:

\begin{equation*}
\begin{aligned}
d \tau &= 4 \left. \delta^1 \wedge \tau \right. + \frac{\hh + \kk}{\A^4} \left. \tau \wedge \sigma \right.
+ \left.\sigma \wedge \zeta \right.
+ \left. \sigma \wedge \ov{\zeta} \right., \\
d \sigma &= 3 \left. \delta^1 \wedge
 \sigma \right. + \frac{\hh + \kk}{\A^4} \left. \tau \wedge \rho \right. + \left. \rho \wedge \zeta \right. 
+ \left.  \rho \wedge \ov{\zeta} \right., \\
d \rho &  =  2 \left. \delta^1 \wedge \rho \right. + i \, \frac{\kk}{\A^4} \, \left. \tau \wedge \zeta \right. 
- i \, \frac{\hh}{\A^4} \left. \tau \wedge \ov{\zeta} \right. + i \, \left. \zeta \wedge \ov{\zeta} \right., \\
d \zeta &= \left. \delta^1 \wedge \zeta 
\right. + \left. \delta^2 \wedge \tau \right. + \frac{\hh \left( \hh + \kk \right)}{\A^8} \left. \tau \wedge \sigma \right.
+
\frac{\hh}{\A^4} \left. \sigma \wedge \zeta \right. + \frac{\hh}{\A^4} \left. \sigma \wedge \ov{\zeta} \right., \\  
d \ov{\zeta} &= \left. \delta^1 \wedge \ov{\zeta} \right. + \left. \delta^3 \wedge \tau \right. 
+ \frac{\kk \left( \hh + \kk \right)}{\A^8} \left. \tau \wedge \sigma \right.
+
\frac{\kk}{\A^4} \left. \sigma \wedge \zeta \right. + \frac{\kk}{\A^4} \left. \sigma \wedge \ov{\zeta} \right. 
.\end{aligned}
\end{equation*}
We deduce from these equations that we can perform the normalization:
\begin{equation*}
\hh = \kk = 0.
\end{equation*}
With the $1$-dimensional group $G_5$ of the form:
\begin{equation*}
g := \begin{pmatrix}
{\A^4}  & 0 & 0 & 0 & 0 \\
0  & \A^3 & 0 & 0 & 0 \\
0 & 0 & \A^2 & 0 & 0 \\
0  & 0 & 0 & \A & 0 \\
0  & 0 & 0 & 0 & \A
\end{pmatrix},
\end{equation*}
whose Maurer-Cartan form is given by 
\begin{equation*}
\alpha:= \frac{d \A}{\A},
\end{equation*}
we get the following structure equations:
\begin{equation*}
\begin{aligned}
d \tau &= 4 \left. \alpha \wedge \tau \right.
+ \left.\sigma \wedge \zeta \right.
+ \left. \sigma \wedge \ov{\zeta} \right., \\
d \sigma &= 3 \left. \alpha \wedge \sigma \right. + \left. \rho \wedge \zeta \right. 
+ \left.  \rho \wedge \ov{\zeta} \right., \\
d \rho &  =  2 \left. \alpha \wedge \rho \right. + i \, \left. \zeta \wedge \ov{\zeta} \right., \\
d \zeta &= \left. \alpha \wedge \zeta \right., \\
d \ov{\zeta} &= \left. \alpha \wedge \ov{\zeta} \right..
\end{aligned}
\end{equation*}
No more normalizations are allowed  
at this stage. We thus just perform a prolongation by adjoining the form $\alpha$ to the structure equations, whose
exterior derivative is given by:
\begin{equation*}
d \alpha = 0.
\end{equation*}
This completes the proof of Theorem~\ref{thm:N}
.

\section{Class ${\sf IV_{2}}$}
Class ${\sf IV_{2}}$ is constituted by the $5$-dimensional real hypersurfaces $M^5 \subset \C^3$ which are of
CR-dimension $2$,  whose Levi form is of constant
rank $1$ and which are $2$-nondegenerate, i.e. their Freeman forms
are non-zero. 
The most symmetric manifold of this class is the tube over the future light cone, 
which is defined by the equation:

\begin{equation*}
{\sf LC}: \,\,\,\,\,\,\,\,\,\,\,\,\,\,\,\,\,\,\,\,\,\,\,\,\, 
\left( {\sf Re} \,z_1 \right)^2 -\left( {\sf Re} \,z_2 \right)^2 
- \left( {\sf Re} \,z_3 \right)^2
=
0, \qquad \qquad {\sf Re} \,z_1 > 0
.\end{equation*}
This section is devoted to the determination of the Lie algebra ${\sf aut_{CR}}({\sf LC})$ of infinitesimal CR-automorphisms of ${\sf LC}$.
This has been done before by Kaup and Zaitsev \cite{Kaup-Zaitsev}.
We prove the following result:
\begin{theorem}
\label{thm:LC}
The tube over the future light cone:
\begin{equation*}
{\sf LC}: \,\,\,\,\,\,\,\,\,\,\,\,\,\,\,\,\,\,\,\,\,\,\,\,\, 
\left( {\sf Re} \,z_1 \right)^2 -\left( {\sf Re} \,z_2 \right)^2 
- \left( {\sf Re} \,z_3 \right)^2
=
0, \qquad \qquad {\sf Re} \,z_1 > 0
.\end{equation*}
has a ${\bf 10}$-dimensional Lie algebra of CR-automorphisms. 
A basis for the Maurer-Cartan forms of ${\sf aut_{CR}}({\sf LC})$ is provided
by the $10$ differential $1$-forms  $\rho$, $\kappa$, $\zeta$, $\ov{\kappa}$, $\ov{\zeta}$, $\pi^1$, $\pi^2$, $\ov{\pi^1}$, 
$\ov{\pi^2}$, $\Lambda$, which satisfy the Maurer-Cartan equations:

\begin{equation}
\label{eq:structure}
\begin{aligned}
d \rho & = \pi^1 \wedge \rho + \ov{\pi^1} \wedge \rho + i \, \kappa \wedge \ov{\kappa}, \\
d \kappa & = \pi^1 \wedge \kappa + \pi^2 \wedge \rho + \zeta \wedge \ov{\kappa}, \\
d \zeta & = i \, \pi^2 \wedge \kappa + \pi^1 \wedge \zeta - \ov{\pi^1} \wedge \zeta,  \\
d \ov{\kappa} & = \ov{\pi^1} \wedge \ov{\kappa} + \ov{\pi^2} \wedge \rho - \kappa \wedge \ov{\zeta}, \\
d \ov{\zeta} & = - i \, \ov{\pi^2} \wedge \ov{\kappa} + \ov{\pi^1} \wedge \ov{\zeta} - {\pi^1} \wedge \ov{\zeta},  \\
d \pi^1 & = \Lambda \wedge \rho + i  \, \kappa \wedge \ov{\pi^2} + \zeta \wedge \ov{\zeta}, \\
d \pi^2 & =   \Lambda \wedge \kappa +  \zeta  \wedge \ov{\pi^2} + \pi^2 \wedge \ov{\pi^1},        \\
d \ov{\pi^1} & = \Lambda \wedge \rho - i  \, \ov{\kappa} \wedge \pi^2 - \zeta \wedge \ov{\zeta}, \\
d \ov{\pi^2} & =   \Lambda \wedge \ov{\kappa} +  \ov{\zeta}  \wedge \pi^2 - \pi^1 \wedge \ov{\pi^2}, \\   
d \Lambda & =  - \pi^1 \wedge \Lambda + i\, \pi^2 \wedge \ov{\pi^2} - \pi^1 \wedge \ov{\Lambda}.
\end{aligned}
\end{equation}
\end{theorem}

\subsection{Geometric set-up}
In order to motivate our subsequent notations, it is convenient to introduce some general results on CR-manifolds
belonging to class ${\sf IV_{2}}$, for which we refer to \cite{pocchiola} for a proof.

Let $M \subset \C^3$ be a smooth hypersurface locally represented 
as a graph over the $5$-dimensional real hyperplane
$\C_{ z_1} \times \C_{ z_2} \times \R_v$:
\[
u
=
F\big(z_1,z_2,\overline{z_1},\overline{z_2},v\big),
\]
where $F$ is a local smooth function depending
on $5$ arguments.
We assume that $M$ is a CR-submanifold of CR dimension $2$ 
which is $2$-nondegenerate and whose Levi form is of constant rank $1$.
The two vector fields $\LL$ and $\mathcal{L}_2$ defined by:

\begin{align*}
\mathcal{L}_j &
=
\frac{\partial}{\partial z_j}
+
A^j\,
\frac{\partial}{\partial v},
&&  A^j := - i \, \frac{F_{z_j}}{1 + i \, F_v},
&&& j= 1,2,
\end{align*}
constitute a basis of $T^{1,0}_pM$ at each point $p$ of $M$ and thus
provide an identification of $T^{1,0}_pM$ with $\C^2$ at each
point. Moreover, the real $1$-form $\sigma$ defined by:
\begin{equation*}
\sigma := dv - A^1 \, dz_1 - A^2 \, dz_2 - \ov{A^1} \, d\ov{z_1} 
- \ov{A^2} \, d \ov{z_2},
\end{equation*}
satisfies
\begin{equation*}
\{ \sigma =0 \} = T^{1,0}M \oplus T^{0,1}M,
\end{equation*}
and thus provides an identication of the projection
\begin{equation*}
\C \otimes T_pM \longrightarrow \left. \C \otimes T_pM \right.\big/ 
\left( T^{1,0}_pM \oplus T^{0,1}_pM \right)
\end{equation*}
with the map $\sigma_p$: \, $\C \otimes T_pM \longrightarrow \C$.
With these two identifications, the Levi form $LF$ can be viewed at
each point $p$ as a skew hermitian form on $\C^2$ represented by the matrix:
\begin{equation*}
LF= 
\begin{pmatrix}
\sigma_p \left( i \, \big[ \LL , \Lb \big] \right)
& \sigma_p \left( i \,\big[ \mathcal{L}_2, {\Lb} \big] \right)\\
\sigma_p \left( i \,\big[ \LL , \ov{\mathcal{L}_2} \big] \right) 
& \sigma_p \left( i \,\big[ \mathcal{L}_2, \ov{\mathcal{L}_2} \big] \right) \\
\end{pmatrix}
.\end{equation*}
The fact that $LF$ is supposed to be of constant rank $1$ ensures the
existence of a certain function $k$ such that the vector field 
\begin{equation*}
\KK := k \,
\LL + \mathcal{L}_2
\end{equation*}
lies in the kernel of $LF$. Here are the expressions 
of $\KK$ and $k$ in terms of the graphing function $F$:
\begin{equation*}
\KK = k \, \partial_{z_1} +  \partial_{z_2}  - \frac{i}{1 + i \, F_v}
  \left( k \, F_{z_1} + F_{z_2} \right) \partial_v,
\end{equation*}
{\tiny
\begin{equation*}
k = - \frac{F_{z_2, \ov{z_1}} + F_{z_2, \ov{z_1}} 
\,  F_{v}^2 - i \, F_{\ov{z_1}} \,  F_{z_2, v} - F_{\ov{z_1}} 
\,  F_{v} \, F_{v, z_2} + i \, F_{z_2} \,  F_{\ov{z_1}} \, 
F_{v,v} - F_{z_2} \,  F_v \,  F_{v, \ov{z_1}}}{
F_{z_1, \ov{z_1}} + F_{z_1, \ov{z_1}} \,  F_v^2 - i \, F_{\ov{z_1}}
 \,  F_{z_1, v} - F_{\ov{z_1}} \,  F_v \,  F_{ z_1, v} + i \, F_{z_1} 
\, F_{\ov{z_1}, v} 
+ F_{z_1} \, F_{\ov{z_1}} \,  F_{v, v} -F_{z_1} \, F_v  \, F_{v,\ov{z_1}}}
.\end{equation*}}
From the above construction, the four vector fields $\LL$, $\KK$, $\Lb$,
$\Kb$ constitute a basis of $T^{1,0}_pM \oplus T^{0,1}_pM$ at each
point $p$ of $M$.  It turns out that the vector field $\T$ defined by:
\begin{equation*}
\T:= i \, \big[ \LL ,\Lb \big]
\end{equation*}
is linearly independant from $\LL$,
$\KK$, $\Lb$, $\Kb$.

It is well known (see \cite{Fels-Kaup-2007, Merker-2003}) that the tube over the future light cone is
locally biholomorphic to
the graphed hypersurface:
\begin{equation*}
u = \frac{z_1 \ov{z_1} + \frac{1}{2} z_1^2 \ov{z_2} + \frac{1}{2} \ov{z_1^2} z_2}{1 - z_2 \ov{z_2}}. 
\end{equation*}
The five vector fields  $\LL$, $\KK$, $\Lb$,
$\Kb$ and $\T$, wich constitute a local frame on $\LC$, have thus the following expressions:

\begin{equation*}
\LL:=
\frac{\partial}{\partial z_1}
-
i \frac{\ov{z_1} + z_1   \ov{z_2}}{1 - z_2  \ov{z_2}} 
\frac{\partial}{\partial v}
, \end{equation*} 

\begin{equation*}
\KK:=  - \frac{\ov{z_1} + z_1   \ov{z_2}}{1 - z_2  
\ov{z_2}} \, \frac{\partial}{\partial z_1} + \frac{\partial}{\partial z_2} + 
\frac{i} {2}  \frac{\ov{z_1}^2 +2 z_1   \ov{z_1}   \ov{z_2} + z_1^2   \ov{z_2}^2}{\left(1 - z_2  \ov{z_2}\right)^2} \,
\frac{\partial}{\partial v},
\end{equation*}
and 

\begin{equation*}
\T := -\frac{2}{1 - z_2 \ov{z_2}} \,  \frac{\partial}{\partial v}.
\end{equation*}
Moreover the function $k$ is given by

\begin{equation*}
k:= - \frac{\ov{z_1} + z_1   \ov{z_2}}{1 - z_2  \ov{z_2}}
.\end{equation*}

Let $(\rho_0, \kappa_0, \zeta_0, \ov{\kappa_0}, \ov{\zeta}_0)$ 
be the dual coframe of $(\T, \LL, \KK, \Lb, \Kb)$. We have:
\begin{dgroup*}
\begin{dmath*}
\rho_0 =- \frac{i}{2} \, \left( \ov{z_1} + z_1 \ov{z_2} \right) \, dz_1 
- \frac{i} {4}  \frac{\ov{z_1}^2 +2 z_1   \ov{z_1}   \ov{z_2} + z_1^2   \ov{z_2}^2}{1 - z_2  \ov{z_2}} \, dz_2
+ \frac{i}{2} \, \left( z_1 + \ov{z_1} z_2 \right) \, d \ov{z_1} 
 + \frac{i} {4}  \frac{z_1^2 +2 z_1 z_2   \ov{z_1}  + \ov{z_1}^2  z_2^2}{1 - z_2  \ov{z_2}} \, dz_2 + \frac{1}{2} \, \left(-1 + z_2 \ov{ z_2} \right) \, dv 
\end{dmath*},
\begin{dmath*}
\kappa_0 = dz_1 + \frac{\ov{z_1} + z_1 \ov{z_2}}{1 - z_2 \ov{z_2}} dz_2
\end{dmath*},
\begin{dmath*}
\zeta_0 = dz_2
\end{dmath*},
\begin{dmath*}
\ov{\kappa_0} =  d \ov{z_1} + \frac{z_1 +\ov{ z_1} z_2}{1 - z_2 \ov{z_2}} d \ov{z_2}
\end{dmath*},
\begin{dmath*}
\ov{\zeta_0} = d \ov{z_2}
\end{dmath*}.
\end{dgroup*}

A direct computation gives the structure equations enjoyed by the coframe 
$(\rho_0, \kappa_0, \zeta_0, \ov{\kappa_0}, \ov{\zeta}_0)$:

\begin{equation}
\label{eq:StEq}
\begin{aligned}
d \rho_{0} & =   \frac{\ov{z_2}}{1-z_2 \ov{z_2}} \left. \rho_0 \wedge \zeta_0 \right.
+ \frac{z_2}{1-z_2 \ov{z_2}}  \left. \rho_0 \wedge \ov{\zeta_0} \right. + i \left.  \kappa_0 \wedge \ov{\kappa_0} \right., \\
d \kappa_{0} & =  \frac{\ov{z_2}}{1-z_2 \ov{z_2}} \, \left. \kappa_0 \wedge \zeta_0 \right. 
- \frac{1}{1- z_2 \ov{z_2}} \, \left. \zeta_0 \wedge \ov{\kappa}_0 \right. 
, \\
d \zeta_0 & = 0, \\
d \ov{\kappa}_0 & =
\frac{1}{1- z_2 \ov{z_2}} \, \left. \kappa_0 \wedge 
 \ov{\zeta_0} \right. + \frac{z_2}{1-z_2 \ov{z_2}} \, \left. \ov{\kappa_0} \wedge \ov{\zeta_0} \right. 
,\\
d \ov{ \zeta_{0}} & = 0.
\end{aligned}
\end{equation}

The matrix Lie group which encodes the equivalence problem for $\LC$ is the  $10$ dimensional Lie group $G_1$ whose elements are of the form:
\begin{equation*}
g := \begin{pmatrix}
{\sf c}\overline{\sf c} & 0 & 0 & 0 & 0 \\
{\sf b} & {\sf c} & 0 & 0 & 0 \\
{\sf d} & {\sf e} & {\sf f} & 0 & 0 \\
\overline{\sf b} & 0 & 0 & \overline{\sf c} & 0 \\
\overline{\sf d} & 0 & 0 & \overline{\sf e} & \overline{\sf f}
\end{pmatrix},
\end{equation*}
where $\cc$ and $\ff$ are non-zero complex numbers whereas $\bb$, $\dd$ and $\ee$ are arbitrary complex numbers (see \cite{pocchiola, MPS}).
We introduce the $5$ new one-forms
$\rho$, $\kappa$, $\zeta$, $\overline{ \kappa}$, $\overline{ \zeta}$ by the relation:

\begin{equation*}
\begin{pmatrix}
\rho
\\
\kappa
\\
\zeta
\\
\overline{\kappa}
\\
\overline{\zeta}
\end{pmatrix}
:=
g \cdot 
\begin{pmatrix}
\rho_0
\\
\kappa_0
\\
\zeta_0
\\
\overline{\kappa_0}
\\
\overline{\zeta_0}
\end{pmatrix}
,
\end{equation*}
which we abbreviate as:
\begin{equation*}
\omega := g \cdot \omega_0.
\end{equation*}

The coframes $\omega$ define a $G_1$ structure $P^1$ on $\LC$. The rest of this section is devoted to reduce $P^1$ 
to an absolute parallelism on $\LC$ through Cartan 
equivalence method.

\subsection{Normalization of $\ff$}
The Maurer Cartan forms of $G_1$ are the following:
\begin{equation*}
\begin{aligned}
\alpha^1  &: = \frac{ d \cc}{\cc}, \\
\alpha^2 &:= \frac{ d \bb}{\cc \cb} - \frac{\bb \, d \cc}{\cc^2}{\cb}, \\
\alpha^3  &:= \frac{ d \dd}{\cc \cb} - \frac{\bb \, 
 d \ee}{\cc^2 \cb} + \frac{\left( -\dd \cc + \ee \bb \right) d {\sf f}}{\cc^2 \cb {\sf f}}, \\
\alpha^4 & := \frac{d \ee}{\cc} - \frac{\ee \, d {\sf f}}{\cc \sf f}, \\
\alpha^5 &:= \frac{d {\sf f}}{\sf f}.
\end{aligned}
\end{equation*}

The structure equations read as:
\begin{dgroup*}
\begin{dmath*}
d \rho  =
\alpha^1 \wedge \rho + \overline{\alpha^{1}} \wedge \rho \\
+
T^{\rho}_{\rho \kappa} \, \rho \wedge \kappa
+
T^{\rho}_{\rho \zeta} \, \rho \wedge \zeta
+
T^{\rho}_{\rho \overline{\kappa}} \, \rho \wedge \overline{\kappa}
+
T^{\rho}_{\rho \overline{\zeta}} \, \rho \wedge \overline{\zeta}
+
i \, \kappa \wedge \overline{\kappa}
,\end{dmath*}
\begin{dmath*}
d \kappa
=
\alpha^{1} \wedge \kappa + \alpha^{2} \wedge \rho \\
+
T^{\kappa}_{\rho \kappa} \, \rho \wedge \kappa
+
T^{\kappa}_{\rho \zeta} \, \rho \wedge \zeta
+
T^{\kappa}_{\rho \overline{\kappa}} \, \rho \wedge \overline{\kappa}
+
T^{\kappa}_{\rho \overline{\zeta}} \, \rho \wedge \overline{\zeta}
+
T^{\kappa}_{\kappa \zeta} \, \kappa \wedge \zeta 
+
T^{\kappa}_{\kappa \overline{\kappa}} \, \kappa \wedge \overline{\kappa}
+
T^{\kappa}_{\zeta \overline{\kappa}} \, \zeta \wedge \overline{\kappa},
\end{dmath*}
\begin{dmath*}
d \zeta
=
\alpha^{3} \wedge \rho + \alpha^{4} \wedge \kappa + \alpha^{5} \wedge \zeta \\
+
T^{\zeta}_{\rho \kappa} \, \rho \wedge \kappa
+
T^{\zeta}_{\rho \zeta} \, \rho \wedge \zeta
+
T^{\zeta}_{\rho \overline{\kappa}} \, \rho \wedge \overline{\kappa}
+
T^{\zeta}_{\rho \overline{\zeta}} \,\left. \rho \wedge \overline{\zeta} \right.
+
T^{\zeta}_{\kappa \zeta} \, \kappa \wedge \zeta 
+
T^{\zeta}_{\kappa \overline{\kappa}} \, \kappa \wedge \overline{\kappa}
+
T^{\zeta}_{\zeta \overline{\kappa}} \, \zeta \wedge \overline{\kappa},
\end{dmath*}
\end{dgroup*}
where the expressions of the torsion coefficients $T^{\smallbullet}_{\smallbullet \smallbullet}$ are given in the appendix.

We now proceed with the absorption step of Cartan's method. 
We introduce the modified Maurer-Cartan forms $\widetilde{\alpha}^i$, 
which are a related to the $1$-forms $\alpha^i$ by the relations:
\begin{equation*}
\widetilde{\alpha}^i := \alpha^i -  x_{\rho}^i \, \rho \, - 
x_{\kappa}^i \, \kappa - x_{\zeta}^i \, \zeta \, - \, x_{\overline{\kappa}}^i \, \overline{\kappa} \, - \, 
x_{\overline{\zeta}}^i \, \overline{\zeta},
\end{equation*}
where $x^1$, $x^2$, $x^3$, $x^4$ and $x^5$ are arbitrary complex-valued functions.
The previously written structure equations take the new form:
\begin{multline*}
 d \rho = \widetilde{\alpha}^1 \wedge \rho + \overline{\widetilde{\alpha}^1} \wedge \rho \\
+ \left( T_{\rho \kappa}^{\rho}- x_\kappa^1 -
 x_{\overline{\kappa}}^1 \right)  \rho \wedge \kappa \, + \, \left(T_{\rho \zeta}^{\rho} - x_{\kappa}^1 - 
\overline{  x_{\overline{\zeta}}^1} \right) 
 \rho \wedge \zeta  \\  \,  + \, \left(T_{\rho \overline{\kappa}}^{\rho} - x_{\overline{\kappa}}^1-    
 \overline{ x_{\kappa}^1} \right)  \rho \wedge \overline{\kappa} 
 \, + \, \left( T_{\rho \overline{\zeta}}^{\rho} - x_{\zeta}^1 - x_{\overline{\zeta}}^1 \right)
 \rho \wedge \overline{\zeta} \\  + i \, \kappa \wedge \overline{\kappa},
\end{multline*}
\begin{multline*}
d \kappa = \widetilde{\alpha}^{1} \wedge \kappa + \widetilde{\alpha}^{2} \wedge \rho \\
+ \, \left(T_{\rho \kappa}^{\kappa} - x_\kappa^2 + x_\rho^1 
    \right)  \rho \wedge \kappa \, + \, \left(T_{\rho \zeta}^{\kappa} - x_{\kappa}^2 \right)
\, \rho \wedge \zeta \\ 
+ \, \left( T_{\rho \overline{\kappa}}^{\kappa} - x_{\overline{\kappa}}^2 \right)  \rho \wedge \overline{\kappa}  + 
\left( T^{\kappa}_{\rho \overline{\zeta}} - x_{\overline{\zeta}}^2 \right)  \rho \wedge \overline{\zeta} \\
+ \, \left( T_{\kappa \zeta}^{\kappa} + x_{\zeta}^1 
\right)\, \kappa \wedge \zeta \, 
 + \, \left( T_{\kappa \overline{\kappa}}^{\kappa} - x_{\overline{\kappa}}^1 \right)  \kappa \wedge \overline{\kappa} \\
+ \, T_{\zeta \overline{\kappa}}^{\kappa} \,  \zeta \wedge \overline{\kappa} \, + \, \left( T_{\kappa \overline{\zeta}}^1 
- x_{\kappa \overline{\zeta} }^1 \right) \, \kappa \wedge \zeta,
\end{multline*}
\begin{multline*}
 d \zeta = \widetilde{\alpha}^3 \wedge \rho + \widetilde{\alpha}^4 \wedge \kappa + \widetilde{\alpha}^5 \wedge \zeta  \\
+ \left( T_{\rho \kappa}^{\zeta} 
 - x_{\kappa}^3 + x_{\rho}^4 \right)  \rho \wedge \kappa  +  \left( T_{\rho \zeta}^{\zeta} - x_{\zeta}^3 + x_{\rho}^5 
\right)  \rho \wedge \zeta \\
 + \left(T_{\rho \kappa}^{\zeta} - x_{\overline{\kappa}}^3  \right) \left. \rho \wedge \overline{\kappa} \right.  + 
\left( T_{\rho \overline{\zeta}}^{\zeta} -x_{\overline{\zeta}}^3 \right)  \rho \wedge \overline{\zeta}  \\
 + \, \left( T_{\kappa \overline{\kappa}}^{\zeta} - 
x_{\overline{\kappa}}^4 \right)  \kappa \wedge \overline{\kappa}  + 
\left( T_{\zeta 
\overline{\kappa}}^{\zeta} - x_{\overline{\kappa}}^5 \right)  \zeta \wedge \overline{\kappa} \\
+ \left( x_{\kappa}^5 - x_{\zeta}^4 \right)  \kappa \wedge \zeta - x_
{\overline{\kappa}}^4  \,
 \kappa \wedge \overline{\kappa} \\
\, + \, \left(x_{\overline{\kappa}}^5 - x_{\overline{\zeta}}^4 \right)  \overline{\kappa} \wedge \zeta - x_{\overline{\zeta}}^5 
 \, \zeta \wedge \overline{\zeta}.
\end{multline*}

We then choose  $x^1$, $x^2$, $x^3$, $x^4$ and $x^5$ in a way that eliminates as many torsion coefficients as possible.
We easily see that the only coefficient 
which can not be absorbed is the one in front of $\zeta \wedge \ov{\kappa}$ in $d \kappa$, because it does not depend 
on the $x^i$'s.
We choose the normalization 
\begin{equation*}
T_{\zeta \overline{\kappa}}^{\kappa} = 1,
\end{equation*}
which yields to :
\begin{equation*}
\sf{f} = - \frac{c}{\overline{c}} \, \frac{1}{1- z_2 \ov{z_2}}.
\end{equation*}
We notice that the absorbed structure equations take the form:
\begin{align*}
d \rho & = \widetilde{\alpha}^1 \wedge 
\rho + \overline{\widetilde{\alpha}^1} \wedge \rho + i \, \kappa \wedge \overline{\kappa}, \\
d \kappa & =  \widetilde{\alpha}^{1} \wedge \kappa + \widetilde{\alpha}^{2} \wedge \rho  + \zeta \wedge \overline{\kappa}, \\
d \zeta &  = \widetilde{\alpha}^3 \wedge \rho + \widetilde{\alpha}^4 \wedge \kappa + \widetilde{\alpha}^5 \wedge \zeta. 
\end{align*}

The normalization of ${\sf f}$ 
gives the new relation :
\begin{equation*}
\begin{pmatrix}
\rho \\
\kappa \\
\zeta \\
\overline{\kappa}\\
\overline{\zeta} 
\end{pmatrix}
=
\begin{pmatrix}
{\sf c} \overline{\sf c} & 0 &0 &0 &0 \\
{\sf b} & {\sf c} & 0 & 0 & 0 \\
{\sf d} & {\sf e} &  \frac{c}{\overline{c}} \, \frac{1}{-1+ z_2 \ov{z_2}} & 0 & 0 \\
{\sf \overline{b}} & 0 & 0 & {\sf \overline{c}} & 0 \\
0 & 0 & {\sf \overline{d}} & {\sf \overline{e}} & \frac{\cb}{\cc} \, \frac{1}{-1+ z_2 \ov{z_2}}
\end{pmatrix}
\cdot
\begin{pmatrix}
\rho_0 \\
\kappa_0 \\
\zeta_0 \\
\overline{\kappa}_0 \\
\overline{\zeta}_0
\end{pmatrix}
.\end{equation*}
We thus introduce the new one-form 
\begin{equation*}
\hat{\zeta}_{0} = - \frac{1}{1- z_2 \ov{z_2}} \cdot \zeta_0
,\end{equation*}
such that the previous equation rewrites :
\begin{equation*}
\begin{pmatrix}
\rho \\
\kappa \\
\zeta \\
\overline{\kappa}\\
\overline{\zeta} 
\end{pmatrix}
=
\begin{pmatrix}
{\sf c} \overline{\sf c} & 0 &0 &0 &0 \\
{\sf b} & {\sf c} & 0 & 0 & 0 \\
{\sf d} & {\sf e} & \frac{\sf c}{\overline{\sf c}} & 0 & 0 \\
{\sf \overline{b}} & 0 & 0 & {\sf \overline{c}} & 0 \\
0 & 0 & {\sf \overline{d}} & {\sf \overline{e}} & \frac{\sf c}{\sf \overline{c}}
\end{pmatrix}
\cdot
\begin{pmatrix}
\rho_0 \\
\kappa_0 \\
\hat{\zeta_0} \\
\overline{\kappa}_0 \\
\overline{\hat{\zeta}}_0
\end{pmatrix}.
\end{equation*}

We have reduced the
 $G_{1}$ equivalence problem to a $G_2$ equivalence problem, where $G_2$ is the $8$ dimensional real matrix 
Lie group whose elements are of the form
\begin{equation*}
g = \begin{pmatrix}
{\sf c} \overline{\sf c} & 0 &0 &0 &0 \\
{\sf b} & {\sf c} & 0 & 0 & 0 \\
{\sf d} & {\sf e} & \frac{\sf c}{\overline{\sf c}} & 0 & 0 \\
{\sf \overline{b}} & 0 & 0 & {\sf \overline{c}} & 0 \\
0 & 0 & {\sf \overline{d}} & {\sf \overline{e}} & \frac{\sf c}{\sf \overline{c}}
\end{pmatrix}.
\end{equation*}
We determine the new structure equations enjoyed by the base coframe 
$(\rho_0, \kappa_0, \hat{\zeta_0}, \kappa_0, \overline{\hat{\zeta_0}})$.
We get :
\begin{dgroup*}
\begin{dmath*}
d \rho_{0} =- \ov{z_2} \, \rho_{0} \wedge \hat{\zeta}_0 - z_2 \, \rho_0 \wedge \overline{\hat{\zeta}_{0}} +
 i\, \kappa_0 \wedge \overline{\kappa_0},
\end{dmath*}
\begin{dmath*}
d\kappa_{0}= -\ov{z_2} \,  \kappa_0 \wedge \hat{\zeta_0}  + \hat{\zeta}_{0} \wedge \ov{\kappa_{0}} 
,\end{dmath*}

\begin{dmath*}
d \hat{\zeta}_0 = z_2 \,  \hat{\zeta}_0 \wedge \ov{\hat{\zeta}_0}
\end{dmath*}.
\end{dgroup*}

\subsection{Normalization of $\bb$}
The Maurer forms of the $G_{2}$ 
are given by the independant entries of the matrix $d g \cdot g^{-1}$. 
We have:
\begin{equation*}
 dg \cdot g^{-1}
= 
\begin{pmatrix}
\beta^1 + \ov{\beta^1} & 0 &0 &0 &0 \\
\beta^2 & \beta^1 & 0 & 0& 0 \\
\beta^3 & \beta^4 & \beta^1 - \ov{\beta^1} & 0 &0 \\
\ov{\beta^2} & 0 & 0 & \ov{\beta^1} & 0 \\
\ov{\beta^3} & 0 & 0 &  \ov{\beta^4} & - \beta^1 + \ov{\beta^1}
\end{pmatrix}
,\end{equation*}
where the forms $\beta^1$, $\beta^2$, $\beta^3$ and $\beta^4$ are defined by
\begin{gather*}
\beta^1  : = \frac{d {\sf c}}{\sf c}, \\
\beta^2 := \frac{d \bb}{\cc \cb}- \frac{\bb d \cc}{\cc^2  \cb}, \\
\beta^3 := \frac{\left( - \dd \cc + \ee \bb \right) d \cc}{\cc^3 \cb} 
- \frac{\left(- \dd \cc + \ee \bb \right) d \cb}{\cc^2 \cb^2} + \frac{d \dd}{\cc \cb} -
 \frac{\bb d \ee}{\cc^2 \cb}, \\
\beta^4 := - \frac{\ee d \cc}{\cc^2} + \frac{\ee d \cb}{\cb \cc} + \frac{d \ee}{\cc }.
\end{gather*}
Using formula (\ref{eq:str}), we get the structure 
equations for the lifted coframe $(\rho, \kappa, \zeta, \ov{\kappa}, \ov{\zeta})$ from those of the base coframe
$(\rho_0, \kappa_0, \hat{\zeta}_0, \ov{\kappa_0}, \ov{\hat{\zeta}_0})$:

\begin{multline*}
d \rho
=
\beta^{1} \wedge \rho + \overline{\beta^{1}} \wedge \rho \\
+
U^{\rho}_{\rho \kappa} \, \rho \wedge \kappa
+
U^{\rho}_{\rho \zeta} \, \rho \wedge \zeta
+
U^{\rho}_{\rho \overline{\kappa}} \, \rho \wedge \overline{\kappa} \\
+
U^{\rho}_{\rho \overline{\zeta}} \, \rho \wedge \overline{\zeta}
+
i \, \kappa \wedge \overline{\kappa}
,\end{multline*}

\begin{multline*}
d \kappa
=
\beta^{1} \wedge \kappa + \beta^{2} \wedge \rho \\
+
U^{\kappa}_{\rho \kappa} \, \rho \wedge \kappa
+
U^{\kappa}_{\rho \zeta} \, \rho \wedge \zeta 
+
U^{\kappa}_{\rho \overline{\kappa}} \, \rho \wedge \overline{\kappa} 
+
U^{\kappa}_{\rho \overline{\zeta}} \, \rho \wedge \overline{\zeta} \\
+
U^{\kappa}_{\kappa \zeta} \, \kappa \wedge \zeta 
+
U^{\kappa}_{\kappa \overline{\kappa}} \, \kappa \wedge \overline{\kappa}
+
\zeta \wedge \overline{\kappa}
,\end{multline*}

\begin{multline*}
d \zeta
=
\beta^{3} \wedge \rho + \beta^{4} \wedge \kappa + \beta^{1} \wedge \zeta - \overline{\beta^{1}} \wedge \zeta \\
+
U^{\zeta}_{\rho \kappa} \, \rho \wedge \kappa
+
U^{\zeta}_{\rho \zeta} \, \rho \wedge \zeta
+
U^{\zeta}_{\rho \overline{\kappa}} \, \rho \wedge \overline{\kappa} \\
+
U^{\zeta}_{\rho \overline{\zeta}} \, \rho \wedge \overline{\zeta} 
+
U^{\zeta}_{\kappa \zeta} \, \kappa \wedge \zeta
+
U^{\zeta}_{\kappa \overline{\kappa}} \, \kappa \wedge \overline{\kappa} \\
+
U^{\zeta}_{\kappa \overline{\zeta}} \, \kappa \wedge \overline{\zeta} 
+
U^{\zeta}_{\zeta \overline{\kappa}} \, \zeta \wedge \overline{\kappa} 
+
U^{\zeta}_{\zeta \overline{\zeta}} \, \zeta \wedge \overline{\zeta}
.\end{multline*}

We introduce the modified Maurer-Cartan forms $\widetilde{\beta}^i$ which differ from the $\beta^i$ by 
a linear combination of the $1$-forms $\rho$, $\kappa$, $\zeta$, $\overline{\kappa}$, 
$\overline{\zeta}$, i.e. that is :
\begin{equation*}
\widetilde{\beta}^i= \beta^i -
 y_{\rho}^i \, \rho \, -  y_{\kappa}^i \, \kappa - y_{\zeta}^i \, \zeta \, - \, y_{\overline{\kappa}}^i \, \overline{\kappa} \, - \, 
y_{\overline{\zeta}}^i \, \overline{\zeta}.
\end{equation*}
The structure equations rewrite:
\begin{dgroup*}
\begin{dmath*}
d \rho = \widetilde{\beta}^1 \wedge \rho + \ov{\widetilde{\beta}^1} \wedge \rho \\  +
 \left( U_{\rho \kappa}^{\rho}  - y_{\kappa}^1 - \ov{y}^1_{\ov{\kappa}} \right) \rho \wedge \kappa + 
\left( U_{\rho \zeta}^{\rho}  - y_{\zeta}^1 - \ov{y}^1_{\ov{\zeta}} \right) \rho \wedge \zeta +
 \left( U_{\rho \ov{\kappa}}^{\rho}  - y_{\ov{\kappa}}^1 - \ov{y}^1_{\kappa} \right) \rho \wedge \ov{\kappa} +
\left( U_{\rho \zeta}^{\rho}  - y_{\ov{\zeta}}^1 - \ov{y}^1_{\zeta} \right) \rho \wedge \zeta + i\, \kappa \wedge \ov{\kappa},
\end{dmath*}
\begin{dmath*}
d \kappa = \widetilde{\beta}^1 \wedge \kappa + \widetilde{\beta}^2 \wedge \rho \\
 + \left( U_{\rho \kappa}^{\kappa} + y^1_{\rho} - y^2_{\kappa} \right) \rho \wedge \kappa +  
\left( U_{\rho \zeta}^{\kappa}  - y^2_{\zeta} \right) \rho \wedge \zeta +  \left( U_{\rho \ov{\kappa}}^{\kappa}  - y^2_{\ov{\kappa}} \right) \rho \wedge \ov{\kappa} 
\\+
\left( U_{\rho \ov{\zeta}}^{\kappa}  - y^2_{\ov{\zeta}} \right) \rho \wedge \ov{\zeta} +
\left( U_{\kappa \zeta}^{\kappa}  - y^1_{\zeta} \right) \kappa \wedge \zeta +  
\left( U_{\kappa \ov{\kappa}}^{\kappa}  - y^1_{\ov{\kappa}} \right) \kappa \wedge \ov{\kappa} - y^1_{\ov{\zeta}} \, \kappa \wedge \ov{\zeta} + 
\zeta \wedge \ov{\kappa}
,\end{dmath*}
\begin{dmath*}
d \zeta = \widetilde{\beta}^{3} \wedge \rho + \widetilde{\beta}^{4} \wedge \kappa + \widetilde{\beta}^{1} \wedge \zeta - \overline{\widetilde{\beta}^{1}} \wedge \zeta \\
+ \left( U_{\rho \kappa}^{\zeta} - y^3_{\kappa} + y^4_{\rho} \right) \rho \wedge \kappa +  \left( U_{\rho \zeta}^{\zeta} - 
y^3_{\zeta} + y^1_{\rho} - \ov{y}^1_{\rho} \right) \rho \wedge \zeta +  \left( U_{\rho \ov{\kappa}}^{\zeta} - y^3_{\ov{\kappa}} \right) 
\left.\rho \wedge \ov{\kappa}\right.  +   
\left( U_{\kappa \zeta}^{\zeta} - y^4_{\zeta} + y^1_{\kappa}-\ov{y}^1_{\kappa} \right) \kappa \wedge \zeta +  \left( U_{\kappa \ov{ \kappa}}^{\zeta} -
 y^4_{\ov{\kappa}} \right) \kappa \wedge \ov{\kappa} + \left( U_{\kappa \ov{\zeta}}^{\zeta} - y^4_{\ov{\zeta}} \right) \kappa \wedge \ov{\zeta} +  
\left( U_{\zeta \ov{ \kappa}}^{\zeta} - y^1_{\ov{\kappa}} + \ov{y}^1_{\kappa} \right) \zeta \wedge \ov{\kappa} + 
\left( U_{\zeta \ov{\zeta}}^{\zeta} - y^1_{\ov{\zeta}} + \ov{y}^1_{\zeta} \right) \zeta \wedge \ov{\zeta}   
.\end{dmath*}
\end{dgroup*}

We get the following absorbtion equations:
\begin{alignat*}{3}
y_{\kappa}^1 + \ov{y}^1_{\ov{\kappa}} &= U_{\rho \kappa}^{\rho},  
& \qquad \qquad  y_{\zeta}^1 + \ov{y}^1_{\ov{\zeta}} &=  U_{\rho \zeta}^{\rho},
& \qquad \qquad   y_{\ov{\kappa}}^1 + \ov{y}^1_{\kappa} &=  U_{\rho \ov{\kappa}}^{\rho}, \\
y_{\ov{\zeta}}^1 + \ov{y}^1_{\zeta} & = U_{\rho \zeta}^{\rho}, 
& \qquad \qquad - y^1_{\rho} + y^2_{\kappa}& = U_{\rho \kappa}^{\kappa},  & 
\qquad \qquad   y^2_{\zeta} & =  U_{\rho \zeta}^{\kappa}, \\
 y^2_{\ov{\kappa}} &=U_{\rho \ov{\kappa}}^{\kappa},  & \qquad \qquad  
 y^2_{\ov{\zeta}}& = U_{\rho \ov{\zeta}}^{\kappa}, & \qquad \qquad 
 y^1_{\zeta} &= U_{\kappa \zeta}^{\kappa},  \\  
 y^1_{\ov{\kappa}} & =  U_{\kappa \ov{\kappa}}^{\kappa}, & \qquad \qquad
 y^1_{\ov{\zeta}} &=0, & \qquad \qquad
 y^3_{\kappa} - y^4_{\rho} & = U_{\rho \kappa}^{\zeta},  \\ 
y^3_{\zeta} - y^1_{\rho} + \ov{y}^1_{\rho} &= U_{\rho \zeta}^{\zeta}, & \qquad \qquad 
 y^3_{\ov{\kappa}} & = U_{\rho \ov{\kappa}}^{\zeta}, &\qquad \qquad   
 y^4_{\zeta} - y^1_{\kappa} + \ov{y}^1_{\kappa} &= U_{\kappa \zeta}^{\zeta}, \\
 y^4_{\ov{\kappa}} &= U_{\kappa \ov{ \kappa}}^{\zeta}, &\qquad \qquad 
 y^4_{\ov{\zeta}} &= U_{\kappa \ov{\zeta}}^{\zeta}, &\qquad \qquad   
y^1_{\ov{\kappa}} - \ov{y}^1_{\kappa} &= U_{\zeta \ov{ \kappa}}^{\zeta}, \\ 
y^1_{\ov{\zeta}} - \ov{y}^1_{\zeta} &= U_{\zeta \ov{\zeta}}^{\zeta}.    
\end{alignat*} 
Eliminating the $y^{\bullet}_{\bullet}$ among theses
 equations leads to the following relations between the torsion coefficients :
\begin{align*}
U_{\rho \ov{\kappa}}^{\rho} & = \ov{U_{\rho \kappa}^{\rho}}, \\
U_{\rho \ov{\zeta}}^{\rho} & = \ov{U_{\rho \zeta}^{\rho}}, \\
U_{\rho \zeta}^{\rho}&  = U_{\kappa \zeta}^{\kappa},\\
U_{\zeta \ov{ \zeta}}^{\zeta}&  = - U_{\rho\ov{ \zeta}}^{\rho},   \\
2 \, U_{\kappa \ov{\kappa}}^{\kappa} & =  U_{\zeta \ov{\kappa}}^{\zeta} + U_{\rho \ov{\kappa}}^{\rho}. 
\end{align*}
We verify easily that the first four equations
 do not depend on the group coefficients and are already satisfied. However, the last one does depend on the group 
coefficients. It gives us the normalization of $\bb$ as it rewrites :
\begin{equation*}
\label{eq:nb}
\bb = - i \, \cb \ee 
.\end{equation*} 
The absorbed structure equations rewrite:
\begin{align*}
d \rho &= \widetilde{\beta}^1 \wedge \rho + \ov{\widetilde{\beta}^1} \wedge \rho + i\, \kappa \wedge \ov{\kappa}, \\
d \kappa &= \widetilde{\beta}^1 \wedge \kappa + \widetilde{\beta}^2 \wedge \rho
+ 
\zeta \wedge \ov{\kappa}
,\\
d \zeta &= 
\widetilde{\beta}^{3} \wedge \rho + \widetilde{\beta}^{4} \wedge \kappa 
+ \widetilde{\beta}^{1} \wedge \zeta - \overline{\widetilde{\beta}^{1}} \wedge \zeta  +
\left(  U_{\zeta \ov{\kappa}}^{\zeta} + U_{\rho \ov{\kappa}}^{\rho} -2 \, U_{\kappa \ov{\kappa}}^{\kappa} \right)  
\, \zeta \wedge \ov{\kappa}.
\end{align*}

\subsection{Normalization of $\dd$}
We have thus reduced the group $G_2$ to a new group $G_3$, 
whose elements are of the form 
\begin{equation*}
{\sf g} = \begin{pmatrix}
\cc \cb & 0 & 0 & 0 & 0 \\
- i\, \ee \cb & \cc & 0 & 0 &0 \\
\dd & \ee & \frac{\cc}{\cb} & 0 & 0 \\
i\, \eb \cc & 0 & 0 & \cb & 0 \\
\db & 0 & 0 & \eb & \frac{\cb}{\cc} 
\end{pmatrix}
.
\end{equation*}
It is a six-dimensional real Lie group. We compute its Maurer Cartan forms with the usual formula
\[
d {\sf g} \cdot {\sf g}^{-1}
=
\begin{pmatrix}
\gamma^1 + \ov{\gamma}^1 & 0 & 0 & 0 &0 \\
\gamma^2 & \gamma^1 & 0 & 0 & 0 \\
\gamma^3 & i\, \gamma^2 & \gamma^1 - \ov{\gamma}^1 & 0 & 0 \\
\ov{\gamma}^2 & 0 & 0 & \ov{\gamma}^1 & 0 \\
- \gamma^3 & 0 & 0 &- i\, \ov{\gamma}^2 & - \gamma^1 + \ov{\gamma}^1 \\
\end{pmatrix},
\]
where 
\begin{equation*}
\gamma^1 := \frac{d \cc}{\cc},
\end{equation*}
\begin{equation*}
\gamma^2 := i \, \ee \frac{d \cc}{\cc^2} - i\, \frac{\ee \, d \cb}{\cc \cb} - i\, \frac{d \ee}{\cc},
\end{equation*}
and 
\begin{equation*}
\gamma^3 := \left(\frac{\dd \cc
 + i\, \ee^2 \cb }{\cc^2 \cb}\right) \left(\frac{d \cb}{\cb} - \frac{ d \cc }{\cc} \right) + \frac{d \dd}{\cc \cb} + 
i \, \frac{\ee d \ee}{\cc^2}
.
\end{equation*}

As the normalization of 
$\bb$ does not depend on the base variables, the third loop of Cartan's method is straightforward.
We get the following structure equations:

\begin{multline*}
d \rho
=
\gamma^{1} \wedge \rho +  \overline{\gamma^{1}} \wedge \rho \\
+
V^{\rho}_{\rho \kappa} \, \rho \wedge \kappa
+
V^{\rho}_{\rho \zeta} \, \rho \wedge \zeta 
+
V^{\rho}_{\rho \overline{\kappa}} \, \rho \wedge \overline{\kappa} \\ 
+
V^{\rho}_{\rho \overline{\zeta}} \, \rho \wedge \overline{\zeta} 
+
i \, \kappa \wedge \overline{\kappa}
,\end{multline*}

\begin{multline*}
d \kappa
=
\gamma^{1} \wedge \kappa + \gamma^{2} \wedge \rho \\
+
V^{\kappa}_{\rho \kappa} \, \rho \wedge \kappa
+
V^{\kappa}_{\rho \zeta} \, \rho \wedge \zeta
+
V^{\kappa}_{\rho \overline{\kappa}} \, \rho \wedge \overline{\kappa} \\
+
V^{\kappa}_{\rho \overline{\zeta}} \, \rho \wedge \overline{\zeta}
+
V^{\kappa}_{\kappa \zeta} \, \kappa \wedge \zeta
+
V^{\kappa}_{\kappa \overline{\kappa}} \, \kappa \wedge \overline{\kappa}
+
\zeta \wedge \overline{\kappa}
,\end{multline*}
\begin{multline*}
d \zeta
=
\gamma^{3} \wedge \rho + i \, \gamma^{2} \wedge \kappa + \gamma^{1} \wedge 
\zeta - \overline{\gamma^{1}} \wedge \zeta \\
+
V^{\zeta}_{\rho \kappa}  \, \rho \wedge \kappa
+
V^{\zeta}_{\rho \zeta} \, \rho \wedge \zeta
+
V^{\zeta}_{\rho \overline{\kappa}} \, \rho \wedge \overline{\kappa} 
+
V^{\zeta}_{\rho \overline{\zeta}} \, \rho \wedge \overline{\zeta} \\
+
V^{\zeta}_{\kappa \zeta} \, \kappa \wedge \zeta
+
V^{\zeta}_{\kappa \overline{\kappa}} \, \kappa \wedge \overline{\kappa} 
+
V^{\zeta}_{\kappa \overline{\zeta}} \, \kappa \wedge \overline{\zeta} 
+
V^{\zeta}_{\zeta \overline{\kappa}} \, \zeta \wedge \overline{\kappa} \\
+
V^{\zeta}_{\zeta \overline{\zeta}} \, \zeta \wedge \overline{\zeta}
.\end{multline*}

We now start the absorption step.
We introduce:
\[
\widetilde{\gamma}^i := \gamma^i -  z^i_{\rho} \,  \rho - z^i_{\kappa} 
 \,  \kappa - z^i_{\zeta} \,  \zeta - z^i_{\ov{\kappa}} \,  \ov{\kappa} - z^i_{\ov{\zeta}} \, \ov{\zeta}
.\]

The structure equations are modified accordingly:
\begin{dgroup*}
\begin{dmath*}
d \rho = \widetilde{\gamma}^{1} \wedge \rho +  \overline{\widetilde{\gamma}^{1}} \wedge \rho \\
+ \left( V_{\rho \kappa}^{\rho} - z^1_{\kappa}- \ov{z^1_{\ov{\kappa}}} \right) \left. \rho \wedge \kappa \right. + 
\left( V_{\rho \zeta}^{\rho} - z^1_{\zeta} - \ov{z^1_{\ov{\zeta}}} \right) \left. \rho \wedge \zeta \right.
+ \left( V_{\rho \ov{\kappa}} - z^1_{\ov{\kappa}} -\ov{z^1_{\kappa}} \right) \left. \rho \wedge \ov{\kappa} \right.
+ \left( V^\rho_{\rho \ov{\zeta}} - \ov{z^1_{\zeta}} - z^1_{\ov{\zeta}} \right) \left. \rho \wedge \ov{\zeta} \right.
,\end{dmath*} 

\begin{dmath*}
d \kappa = \widetilde{\gamma}^{1} \wedge \kappa + \widetilde{\gamma}^{2} \wedge \rho \\
+ \left( V_{\rho \kappa}^{\kappa}
 - z^2_{\kappa} + z^{\kappa}_{\rho} \right) \left. \rho \wedge \kappa \right. + \left( V_{\rho \zeta}^{\kappa}
 - z^2_{\zeta} \right) 
\left. \rho \wedge \zeta \right. +
 \left( V^{\kappa}_{\rho \ov{\kappa}} - z^2_{\ov{\kappa}}
 \right) \left. \rho \wedge \ov{\kappa} \right.
+ \left( V_{\rho \ov{\zeta}}- z^2_{\ov{\zeta}}
 \right) \left. \rho \wedge \ov{\zeta} \right. + \left( V_{\kappa \zeta}^{\kappa} - z^1_{\zeta} \right) 
\left. \kappa \wedge \zeta \right. + 
\left( V^{\kappa}_{\kappa \ov{\kappa}} 
- z^1_{\ov{\kappa}} \right) \left. \kappa \wedge \ov{\kappa} \right. 
+ \zeta \wedge \ov{\kappa} - z^1_{\ov{\zeta}} \left. 
\kappa \wedge \ov{\zeta} \right.
,\end{dmath*}
\end{dgroup*}
and 
\begin{dmath*}
d \zeta = \widetilde{\gamma}^{3} \wedge \rho + i \, \widetilde{\gamma}^{2} \wedge \kappa + \widetilde{\gamma}^{1} \wedge 
\zeta - \overline{\widetilde{\gamma}^{1}} \wedge \zeta \\
+ \left( V^{\zeta}_{\rho \kappa} - z^3_{\kappa} + i \, z_{\rho}^2 \right) \left. \rho \wedge \kappa \right. 
+ \left( V^{\zeta}_{\rho \zeta} + z^1_{\rho} - z^3_{\zeta}- \ov{\zeta^1_{\rho}} \right) \left. \rho \wedge \zeta \right.
+ \left( V^{\zeta}_{\rho \ov{\kappa}} - z^3_{\ov{\kappa}} \right) \left. \rho \wedge \ov{\kappa} \right. +
\left( V^{\zeta}_{\rho \ov{\zeta}} - z^3_{\ov{\zeta}} \right) \left. \rho \wedge \ov{\zeta} \right. +
\left( V^{\zeta}_{\kappa \ov{\kappa}} - i \, z^2_{\ov{\kappa}} \right) \left. \kappa \wedge \ov{\kappa} \right.
+ \left( V^{\zeta}_{\kappa \ov{\zeta}} - i \, z^2_{\ov{\zeta}} \right) \left. \kappa \wedge \ov{\zeta} \right.
+ \left( V^{\zeta}_{\zeta \ov{\kappa}} - 
z^1_{\ov{\kappa}} + \ov{z^1_{\kappa}} \right) \left. \zeta \wedge \ov{\zeta} \right.
.\end{dmath*}

We thus want to solve the system of linear equations :

\begin{alignat*}{3}
z^1_{\kappa} + \ov{z^1_{\ov{\kappa}}}  &= V^{\rho}_{\rho \kappa},  
& \qquad \qquad z^1_{\ov{\kappa}} + \ov{z^1_{\kappa}}  &= V^{\rho}_{\rho \ov{\kappa}}, 
& \qquad \qquad z^1_{\zeta} + \ov{z^1_{\ov{\zeta}}}  &= V^{\rho}_{\rho \zeta},  \\
\ov{z^1_{\zeta}} + z^1_{\ov{\zeta}} &= V^{\rho}_{\rho \ov{\zeta}}, 
& \qquad \qquad z^2_{\kappa} - z^1_{\rho} & = V^{\kappa}_{\rho \zeta}, 
& \qquad \qquad z^2_{\ov{\kappa}} &= V^{\kappa}_{\rho \ov{\kappa}}, \\
z^2_{\zeta}& = V^{\kappa}_{\rho \zeta}, 
& \qquad \qquad z^2_{\ov{\zeta}} &= V^{\kappa}_{\rho \ov{\zeta}}, 
& \qquad \qquad z^1_{\zeta}& = V^{\kappa}_{\kappa \zeta}, \\
z^1_{\ov{\zeta}} &= 0, 
& \qquad \qquad z^1_{\ov{\kappa}} &= V^{\kappa}_{\kappa \ov{\kappa}}, 
& \qquad \qquad z^3_{\kappa} - i \, z^2_{\rho} &= V^{\zeta}_{\rho \kappa}, \\
- z^1_{\rho} + \ov{z^1_{\rho}} + z^3_{\zeta} &= V^{\zeta}_{\rho \zeta}, 
& \qquad \qquad z^1_{\kappa} - \ov{z^1_{\ov{\kappa}}} - i\, z^2_{\zeta} &= - V^{\zeta}_{\kappa \zeta}, 
& \qquad \qquad i \, z^2_{\ov{\kappa}} &= V^{\zeta}_{\kappa \ov{\kappa}}, \\
z^3_{\ov{\kappa}}& = V^{\zeta}_{\rho \ov{\kappa}}, 
& \qquad \qquad z^3_{\ov{\zeta}} &= V^{\zeta}_{\rho \ov{\zeta}}, 
& \qquad \qquad i \, z^2_{\ov{\zeta}} &= V^{\zeta}_{\kappa \ov{\zeta}}, \\
z^1_{\kappa} - \ov{z^1_{\kappa}} &= V^{\zeta}_{\zeta \ov{\kappa}},
& \qquad \qquad z^1_{\ov{\zeta}} - \ov{z^1_{\zeta}} &= V^{\zeta}_{\zeta \ov{\zeta}}.
\end{alignat*}

This is easily done as:
\begin{equation*}
\left\{
\begin{aligned}
z^1_{\kappa}  & = \ov{V^{\rho}_{\rho \ov{\kappa}}}, \\
z^1_{\ov{\kappa}} & = V^{\kappa}_{\kappa \ov{\kappa}}, \\
z^1_{\zeta} & = V^{\rho}_{\rho \zeta}, \\
z^1_{\ov{\zeta}} & = 0, \\
z^2_{\ov{\kappa}} & = V^{\kappa}_{\rho \ov{\kappa}}, \\
z^2_{\ov{\zeta}} & = V^{\kappa}_{\rho \ov{\zeta}}, \\
z^2_{\zeta} & = V^{\kappa}_{\rho \zeta}, \\
z^3_{\ov{\kappa}} & = V^{\zeta}_{\rho \ov{\kappa}}, \\
z^3_{\ov{\zeta}} & = V^{\zeta}_{\rho \ov{\zeta}}, \\
z^3_{\zeta} & = V^{\zeta}_{\rho \zeta} + z^1_{\rho} - z^1_{\rho}, \\
z^3_{\kappa} & = V^{\zeta}_{\rho \kappa} + i\, z^2_{\rho}, \\
z^2_{\kappa} & = V^{\kappa}_{\rho \zeta} + z^1_{\rho}, 
\end{aligned}
\right.
\end{equation*}
where $z^1_{\rho}$ and $z^2_{\rho}$ may be choosen freely. Eliminating the $z^{\bullet}_{\bullet}$
we get the following additional conditions on the $V^{\bullet}_{\bullet \bullet}$ :
\begin{equation} \label{eq:31}
\left\{
\begin{aligned}
V^{\rho}_{\rho \ov{\kappa}} &= \ov{V^{\rho}_{\rho \kappa}}, \\
V^{\rho}_{\rho \ov{\zeta}}  &= \ov{V^{\rho}_{\rho \zeta}}, \\
V^{\rho}_{\rho \zeta} &= V^{\kappa}_{\kappa \zeta}, \\
i \, V^{\kappa}_{\rho \ov{\zeta}} &= V^{\zeta}_{\kappa \ov{\zeta}}, \\
V^{\rho}_{\rho \zeta} & = - \ov{V^{\zeta}_{\zeta \ov{\zeta}}}, \\
2 \, V^{\kappa}_{\kappa \ov{\kappa}} & = V^{\rho}_{\rho \ov{\kappa}} + V^{\zeta}_{\zeta \ov{\kappa}},
\end{aligned}
\right.
\end{equation}
and
\begin{equation}
\left\{
\begin{aligned}
i \, V^{\kappa}_{\rho \ov{\kappa}} & = V^{\zeta}_{\kappa \ov{\kappa}},  \\
V^{\ov{\zeta}}_{\kappa \ov{\zeta}} +  V^{\zeta}_{\kappa \zeta} & =   i \, V^{\kappa}_{\rho \zeta}.
\end{aligned}
\right.
\end{equation}

We easily verify that the equations (\ref{eq:31}) 
are indeed satisfied. However the remaining two equations are not
and they provide two essential torsion coefficients, namely $i \, V^{\kappa}_{\rho \ov{\kappa}} - V^{\zeta}_{\kappa \ov{\kappa}}$ and 
$V^{\ov{\zeta}}_{\kappa \ov{\zeta}} +  V^{\zeta}_{\kappa \zeta} -   i \, V^{\kappa}_{\rho \zeta}$,
which will give us at least one new normalization of the group coefficients.
Indeed we have 
\begin{equation*}
i \, V^{\kappa}_{\rho \ov{\kappa}} - V^{\zeta}_{\kappa \ov{\kappa}}= - 2i \, \frac{\dd}{\cc \cb} + \frac{\ee^2}{\cc^2}
.\end{equation*}
Setting this expression to $0$, we get the normalization of the parameter $\dd$:
\begin{equation*}
\dd = - i\, \frac{1}{2} \,  \frac{\ee^2 \cb}{ \cc}.
\end{equation*}

\subsection{Prologation of the $G_4$ structure}
We have reduced the previous
 $G_3$-structure to a $G_4$-structure, where $G_4$ is the four dimensional matrix Lie group whose elements are of the form :
\begin{equation*}
\begin{pmatrix}
\cc \cb & 0 & 0 & 0 & 0 \\
-i \, \ee \cb & \cc & 0 & 0 & 0 \\
- \frac{i}{2} \, \frac{\ee^2 \cb}{\cc}  & \ee & \frac{\cc}{\cb} & 0 & 0 \\
i \,  \eb \cc & 0 & 0 &  \cb & 0  \\
 \frac{i}{2} \, \frac{\eb^2 \cc}{\cb}  &0 & 0 &  \eb & \frac{\cb}{\cc}  \\
\end{pmatrix}
.\end{equation*}

The basis for the Maurer-Cartan forms of $G_4$ is provided by the four forms 
\begin{equation*}
\delta^1 := \frac{d \cc}{\cc} \quad, \quad \delta^2 :=  i \, \ee \frac{d \cc}{\cc^2} 
- i\, \frac{\ee \, d \cb}{\cc \cb} - i\, \frac{d \ee}{\cc} \quad ,
\quad \ov{\delta^1} \quad , \quad \ov{\delta^2}.
\end{equation*}

Now we just substitute the parameter $\dd$ by its 
normalization in the structure equations at the third step. 
We have to take into account the fact that $d \dd$ is modified accordingly. Indeed we have:
\begin{equation*}
d \dd = - i \ee \frac{\cb}{\cc} - \frac{i}{2} \, 
\frac{\ee^2 \cb}{\cc} \left( \frac{ d \cb}{\cb} - \frac{d \cc}{\cc} \right) 
.\end{equation*}
The forms $\gamma^1$ and $\gamma^2$ are 
not modified as they do not involve terms in $d \dd$, but this is not the case for $\gamma^3$ which is transformed into:
\begin{align*}
\gamma^3  & = \frac{d \dd }{\cc \cb} +
 i \, \frac{\ee}{\cc^2} - \frac{\dd \, d \cc}{\cc^2 \cb^2}
 - i \, \ee^2 \frac{d \cc}{\cc^3} + \frac{\dd \, d \cb}{\cc \cb^2} + 
i \,\frac{\ee^2 \, d\cb}{\cb \cc^2} \\
          &=  0.
\end{align*}
The expressions of $d \rho$ and $d \kappa$ are thus 
unchanged from the expressions given by the structure equations at the third step, 
except the fact that we shall replace $\dd$ by 
$ - \frac{i}{2} \frac{\ee^2 \cb}{\cc} + i \, \frac{\cc}{\cb} \, H$ 
in the expression of each torsion coefficient $V_{\bullet \bullet}^{\bullet}$, which we rename $W_{\bullet \bullet}^{\bullet}$,
and the fact that the forms $\gamma^1$ and $\gamma^1$ shall be replaced by the forms
$\delta^1$ and $\delta^2$, that is:

\begin{dmath*}
d \rho
 =
\delta^1 \wedge \rho +  \overline{\delta^{1}} \wedge \rho \\
+
W^{\rho}_{\rho \kappa} \, \rho \wedge \kappa
+
W^{\rho}_{\rho \zeta} \, \rho \wedge \zeta 
+
W^{\rho}_{\rho \overline{\kappa}} \, \rho \wedge \overline{\kappa}
+
W^{\rho}_{\rho \overline{\zeta}} \, \rho \wedge \overline{\zeta}
+
i \, \kappa \wedge \overline{\kappa},
\end{dmath*}
and
\begin{dmath*}
d \kappa
 =
\delta^{1} \wedge \kappa + \delta^{2} \wedge \rho \\
+
W^{\kappa}_{\rho \kappa} \, \rho \wedge \kappa
+
W^{\kappa}_{\rho \zeta} \, \rho \wedge \zeta
+
W^{\kappa}_{\rho \overline{\kappa}} \, \left. \rho \wedge \overline{\kappa} \right.
+
W^{\kappa}_{\rho \overline{\zeta}} \,\left.  \rho \wedge \overline{\zeta} \right.
+
W^{\kappa}_{\kappa \zeta} \, \kappa \wedge \zeta
+
W^{\kappa}_{\kappa \overline{\kappa}} \, \kappa \wedge \overline{\kappa}
+
\zeta \wedge \overline{\kappa}.
\end{dmath*}

The expression of $d \zeta$ is obtained in the same way, setting $\gamma_3$ to zero, and renaming $W^{\bullet}_{\bullet \bullet}$ 
the coefficients $V_{\bullet \bullet}^{\bullet}$ in which one performs the substitution 
$d= - i \frac{1}{2} \, \frac{\ee^2 \cb}{\cc}$: 

\begin{dmath*}
d \zeta
=
 i \, \delta_{2} \wedge \kappa + \delta_{1} \wedge 
\zeta - \overline{\delta_{1}} \wedge \zeta \\
+
W^{\zeta}_{\rho \kappa}  \, \rho \wedge \kappa
+
W^{\zeta}_{\rho \zeta} \, \rho \wedge \zeta
+
W^{\zeta}_{\rho \overline{\kappa}} \, \rho \wedge \overline{\kappa} 
+
W^{\zeta}_{\rho \overline{\zeta}} \, \rho \wedge \overline{\zeta} 
+
W^{\zeta}_{\kappa \zeta} \, \left.  \kappa \wedge \zeta \right.
+
W^{\zeta}_{\kappa \overline{\kappa}} \, \kappa \wedge \overline{\kappa}
+
W^{\zeta}_{\kappa \overline{\zeta}} \, \kappa \wedge \overline{\zeta} 
+
W^{\zeta}_{\zeta \overline{\kappa}} \, \zeta \wedge \overline{\kappa}
+
W^{\zeta}_{\zeta \overline{\zeta}} \, \zeta \wedge \overline{\zeta}
.\end{dmath*}

Let us now proceed with the absorption phase. We make the two substitutions: 
\begin{equation*}
\begin{aligned}
& \delta^1 : = \widetilde{\delta}^1 + 
w^1_{\rho} \, \rho + w^1_{\kappa} \, \kappa + w^1_{\zeta} \, \zeta + 
w^1_{\ov{\kappa}} \, \ov{\kappa} + w^1_{\ov{\zeta}} \, \ov{\zeta}, \\
& \delta^2 : =  \widetilde{\delta}^2 + 
w^2_{\rho} \, \rho + w^2_{\kappa} \, \kappa + w^2_{\zeta} \, \zeta + w^2_{\ov{\kappa}} \, \ov{\kappa} + w^2_{\ov{\zeta}} \, \ov{\zeta}, 
\end{aligned}
\end{equation*}
in the previous equations.
We get:
\begin{dmath*}
d \rho = \widetilde{\delta}^1 \wedge \rho +  
\overline{\widetilde{\delta}^{1}} \wedge \rho \\
+ \left( W_{\rho \kappa}^{\rho} - w^1_{\kappa}- \ov{w^1_{\ov{\kappa}}} \right) \left. \rho \wedge \kappa \right. + 
\left( W_{\rho \zeta}^{\rho} - w^1_{\zeta} - \ov{w^1_{\ov{\zeta}}} \right) \left. \rho \wedge \zeta \right.
+ \left( W_{\rho \ov{\kappa}} - w^1_{\ov{\kappa}} -\ov{w^1_{\kappa}} \right) \left. \rho \wedge \ov{\kappa} \right.
+ \left( W^\rho_{\rho \ov{\zeta}} - \ov{w^1_{\zeta}} - w^1_{\ov{\zeta}} \right) \left. \rho \wedge \ov{\zeta} \right.,
\end{dmath*} 

\begin{dmath*}
d \kappa = \widetilde{\delta}^{1} \wedge 
\kappa + \widetilde{\delta}^{2} \wedge \rho \\
+ \left( W_{\rho \kappa}^{\kappa} -
 w^2_{\kappa} + w^{1}_{\rho} \right) 
\left. \rho \wedge \kappa \right. + \left( W_{\rho \zeta}^{\kappa}
 - w^2_{\zeta} \right) \left. \rho \wedge
 \zeta \right. + \left( W^{\kappa}_{\rho \ov{\kappa}} - w^2_{\ov{\kappa}} \right) \left. \rho \wedge \ov{\kappa} \right.
+ \left( W_{\rho \ov{\zeta}}- w^2_{\ov{\zeta}} \right) 
\left. \rho \wedge \ov{\zeta} \right. + \left( W_{\kappa \zeta}^{\kappa} - w^1_{\zeta} \right) 
\left. \kappa \wedge \zeta \right. +
 \left( W^{\kappa}_{\kappa \ov{\kappa}} - w^1_{\ov{\kappa}} \right) \left. \kappa \wedge \ov{\kappa} \right. 
+ \zeta \wedge \ov{\kappa} - w^1_{\ov{\zeta}} \left. \kappa \wedge \ov{\zeta} \right.,
\end{dmath*}
and 
\begin{dmath*}
d \zeta = i \, \widetilde{\delta_{2}} \wedge \kappa + \widetilde{\delta_{1}} \wedge 
\zeta - \overline{\widetilde{\delta_{1}}} \wedge \zeta \\
+\left( W^{\zeta}_{\rho \kappa}  + i \, w_{\rho}^2 \right) \left. \rho \wedge \kappa \right. 
+ \left( W^{\zeta}_{\rho \zeta} + w^1_{\rho} - \ov{w^1_{\rho}} \right) \left. \rho \wedge \zeta \right. 
+  W^{\zeta}_{\rho \ov{\kappa}} \left. \rho \wedge \ov{\kappa} \right. +
 W^{\zeta}_{\rho \ov{\zeta}}  \left. \rho \wedge \ov{\zeta} \right. +
\left( W^{\zeta}_{\kappa \ov{\kappa}} - i \, w^2_{\ov{\kappa}} \right) \left. \kappa \wedge \ov{\kappa} \right.
+ \left( W^{\zeta}_{\kappa \ov{\zeta}} - i \, w^2_{\ov{\zeta}} \right) \left. \kappa \wedge \ov{\zeta} \right.
+ \left( W^{\zeta}_{\zeta \ov{\kappa}} -
 w^1_{\ov{\kappa}} + \ov{w^1_{\kappa}} \right) \left. \zeta \wedge \ov{\zeta} \right.
.\end{dmath*}
From the last equation, we immediately 
see that $ W^{\zeta}_{\rho \ov{\kappa}}$ and $W^{\zeta}_{\rho \ov{\zeta}} $ are two new essential torsion coefficients.
We find the remaining ones by solving the set of equations:

\begin{alignat*}{3}
w^1_{\kappa} + \ov{w^1_{\ov{\kappa}}}  &= W^{\rho}_{\rho \kappa}, 
& \qquad \qquad w^1_{\ov{\kappa}} + \ov{w^1_{\kappa}}  &= W^{\rho}_{\rho \ov{\kappa}}, 
& \qquad \qquad w^1_{\zeta} + \ov{w^1_{\ov{\zeta}}}  &= W^{\rho}_{\rho \zeta},  \\
\ov{w^1_{\zeta}} + w^1_{\ov{\zeta}} &= W^{\rho}_{\rho \ov{\zeta}}, 
& \qquad \qquad w^2_{\kappa} - w^1_{\rho} & = W^{\kappa}_{\rho \kappa}, 
& \qquad \qquad w^2_{\ov{\kappa}} &= W^{\kappa}_{\rho \ov{\kappa}}, \\
w^2_{\zeta}& = W^{\kappa}_{\rho \zeta}, 
& \qquad \qquad w^2_{\ov{\zeta}} &= W^{\kappa}_{\rho \ov{\zeta}}, 
& \qquad \qquad w^1_{\zeta}& = W^{\kappa}_{\kappa \zeta}, \\
w^1_{\ov{\zeta}} &= 0, 
& \qquad \qquad w^1_{\ov{\kappa}} &= W^{\kappa}_{\kappa \ov{\kappa}}, 
& \qquad \qquad  - i \, w^2_{\rho} &= W^{\zeta}_{\rho \kappa}, \\
- w^1_{\rho} + \ov{w^1_{\rho}} &= W^{\zeta}_{\rho \zeta}, 
& \qquad \qquad w^1_{\kappa} - \ov{w^1_{\ov{\kappa}}} - i\, w^2_{\zeta} &= - W^{\zeta}_{\kappa \zeta}, 
& \qquad \qquad i \, w^2_{\ov{\kappa}} &= W^{\zeta}_{\kappa \ov{\kappa}}, \\ 
w^1_{\kappa} - \ov{w^1_{\kappa}} &= W^{\zeta}_{\zeta \ov{\kappa}},
& \qquad \qquad i \, w^2_{\ov{\zeta}} &= W^{\zeta}_{\kappa \ov{\zeta}},
& \qquad \qquad w^1_{\ov{\zeta}} - \ov{w^1_{\zeta}} &= W^{\zeta}_{\zeta \ov{\zeta}},
\end{alignat*}
which lead easily as before to:
\begin{equation} \label{absorption}
\left\{
\begin{aligned}
w^1_{\kappa}  & = \ov{W^{\rho}_{\rho \ov{\kappa}}}, \\
w^1_{\ov{\kappa}} & = W^{\kappa}_{\kappa \ov{\kappa}}, \\
w^1_{\zeta} & = W^{\rho}_{\rho \zeta}, \\
w^1_{\ov{\zeta}} & = 0, \\
w^2_{\ov{\kappa}} & = W^{\kappa}_{\rho \ov{\kappa}}, \\
w^2_{\ov{\zeta}} & = W^{\kappa}_{\rho \ov{\zeta}}, \\
w^2_{\zeta} & = W^{\kappa}_{\rho \zeta}, \\
w^2_{\kappa} & = W^{\kappa}_{\rho \kappa} + w^1_{\rho}, \\
w^2_{\rho}& = W^{\zeta}_{\rho \kappa}, \\
-w^1_{\rho} + \ov{w^1_{\rho}} & = W_{\rho, \zeta}^{\zeta}.
\end{aligned}
\right.
\end{equation}
Eliminating the $w^{\bullet}_{\bullet}$ from (\ref{absorption}), 
we get one additionnal condition on the $W^{\bullet}_{\bullet \bullet}$ which has not yet been checked, namely
that $W_{\rho, \zeta}^{\zeta}$ shall be purely imaginary.
The computation of  $W_{\rho, \zeta}^{\zeta}$,  $ W^{\zeta}_{\rho \ov{\kappa}}$ and $W^{\zeta}_{\rho \ov{\zeta}} $ gives:
\begin{equation*}
W_{\rho, \zeta}^{\zeta} = i \, \frac{\ee \eb}{\cc \cb}
 - \frac{i}{2} \, \frac{\ee^2 \cb}{\cc^3} \, \ov{z_2} -  \frac{i}{2} \, \frac{\eb^2 \cc}{\cb^3} \, z_2
,\end{equation*}
\begin{equation*}
 W^{\zeta}_{\rho \ov{\kappa}}=0,
\end{equation*}
and
\begin{equation*}
 W^{\zeta}_{\rho \ov{\zeta}}=0
,\end{equation*}
 which indicates that no further normalizations are 
allowed at this stage and that we must perform a prolongation of the problem.
Let us  introduce the modified Maurer Cartan forms of the group $G_4$, namely :
\begin{equation*}
\left\{
\begin{aligned}
& \hat{\delta}^1 : = \delta^1 -  w^1_{\rho}
 \, \rho - w^1_{\kappa} \, \kappa - w^1_{\zeta} \, \zeta - w^1_{\ov{\kappa}} \, \ov{\kappa} - w^1_{\ov{\zeta}} \, \ov{\zeta}, \\
& \hat{\delta}^2 : =  \delta^2 - w^2_{\rho}
 \, \rho - w^2_{\kappa} \, \kappa - w^2_{\zeta} \, \zeta - w^2_{\ov{\kappa}} \, \ov{\kappa} - w^2_{\ov{\zeta}} \, \ov{\zeta} 
,\end{aligned}
\right.
\end{equation*}
where $w^i_{\rho}$, $w^i_{\kappa}$, $w^i_{\zeta}$, $w^i_{\ov{\kappa}}$, $w^i_{\ov{\zeta}}$, $i = 1 , \,  2 $,  are the solutions 
of the system of equations (\ref{absorption}) corresponding to 
$w^1_{\rho} + \ov{w^1_{\rho}} = 0,$ that is :
\begin{equation}
\label{eq:deltatilde}
\left\{
\begin{aligned}
& \hat{\delta}^1 : = \delta^1 + \frac{1}{2} \, V^{\zeta}_{\rho \zeta} \, \rho 
- \ov{V^{\rho}_{\rho \kappa}} \, \kappa - V^{\rho}_{\rho \zeta} \, \zeta - V^{\kappa}_{\kappa \ov{\kappa}} \, \ov{\kappa}, \\
& \hat{\delta}^2 : =  \delta^2 -
 V^{\zeta}_{\rho \kappa} \, \rho - \left( V^{\kappa}_{\rho \kappa} - \frac{1}{2}  V^{\zeta}_{\rho \zeta} \right) \kappa 
- V^{\kappa}_{\rho \zeta} \, \zeta - V^{\kappa}_{\rho \ov{\kappa}} \, \ov{\kappa} - V^{\kappa}_{\rho \ov{\zeta}} \, \ov{\zeta} 
.\end{aligned}
\right.
\end{equation}
We also introduce the modified
 Maurer Cartan forms which correspond to solutions of the system (\ref{absorption}) when ${\sf Re}(w^1_{\rho})$ 
is not necessarily set to zero, namely :
\begin{equation} 
\label{eq:pi}
\left\{
\begin{aligned}
& \pi^1 : =\hat{ \delta}^1 - \Re(w^1_{\rho}) \, \rho,  \\
& \pi^2 : =  \hat{\delta}^2 - \Re(w^1_{\rho}) \, \kappa.
\end{aligned}
\right.
\end{equation}
Let $P^9$ be the nine dimensional 
$G_4$-structure constituted 
 by the set of all coframes of the form $(\rho, \kappa, \zeta, \ov{\kappa}, \ov{\zeta})$ on $M^5$.
The initial coframe $(\rho_0, \kappa_0, \zeta_0, 
\ov{\kappa}_0, \ov{\zeta}_0)$ gives a natural trivialisation $P^9 \stackrel{p} \rightarrow M^5 \times G_{4}$
which allows us to consider any differential 
form on $M^5$ or $G^4$ as a differential form on $P^9$. If $\omega$ is a differential form on $M^5$ for example, 
we just 
consider $p^*( pr_1^*(\omega))$,  
where $pr_1$ is the projection on the first component $M^5 \times G_4 \stackrel{pr_1} \rightarrow M^5$. 
We still denote this form by $\omega$ in the sequel.
Following \cite{Olver-1995}, we introduce the 
two coframes $(\rho, \kappa, \zeta, \ov{\kappa}, \ov{\zeta}, \delta^1, \delta^2, \ov{\delta^1}, \ov{\delta^2})$ 
and $(\rho, \kappa, \zeta, \ov{\kappa}, 
\ov{\zeta}, \pi^1, \pi^2, \ov{\pi^1}, \ov{\pi^2})$ on $P^9$. Setting $ {\sf t}:= - \Re(w^1_{\rho})$, they relate to each other
by the relation:
\begin{equation*}
\begin{pmatrix}
\rho \\
\kappa \\
\zeta \\
\ov{\kappa} \\
\ov{\zeta} \\
\pi^1 \\
\pi^2 \\
\ov{\pi^1} \\
\ov{\pi^2} \\
\end{pmatrix}
= g_{{\sf t}}
\cdot
\begin{pmatrix}
\rho \\
\kappa \\
\zeta \\
\ov{\kappa} \\
\ov{\zeta} \\
\delta^1 \\
\delta^2 \\
\ov{\delta^1} \\
\ov{\delta^2} \\
\end{pmatrix}
,\end{equation*}
where $g_{\sf t}$ is defined by 
\begin{equation*}
g_{\sf t}:=\begin{pmatrix}
1 & 0 & 0 & 0 & 0 & 0 & 0 & 0 & 0 \\
0 & 1 & 0 & 0 & 0 & 0 & 0 & 0 & 0 \\
0 & 0 & 1 & 0 & 0 & 0 & 0 & 0 & 0 \\
0 & 0 & 0 & 1 & 0 & 0 & 0 & 0 & 0\\
0 & 0 & 0 & 0 & 1 & 0 & 0 & 0 & 0\\
t & 0 & 0 & 0 & 0 & 1 & 0 & 0 & 0 \\
0 & t & 0 & 0 & 0 & 0 & 1 & 0 & 0 \\
t & 0 & 0 & 0 & 0 & 0 & 0 & 1 & 0 \\
0 & 0 & 0 & t & 0 & 0 & 0 & 0 & 1 \\
\end{pmatrix}
.
\end{equation*}
The set 
$\left \{ g_{\sf t}, {\sf t} \in \mathbb{R} \right\}$ defines a one-
dimensional Lie group, whose Maurer Cartan form is given by $d {\sf t}$, 
which we rename $\Lambda$ in the sequel.
We now start the reduction step in the equivalence problem on $P^9$. 
From the definition 
of $\pi^1$ and $\pi^2$ 
as the solutions of the absorption equations $(\ref{absorption})$, the five first structure equations read
as 
\begin{equation}
\label{steq}
\begin{aligned}
d \rho & = \pi^1 \wedge \rho + \ov{\pi^1} \wedge \rho + i \, \kappa \wedge \ov{\kappa}, \\
d \kappa & = \pi^1 \wedge \kappa + \pi^2 \wedge \rho + \zeta \wedge \ov{\kappa}, \\
d \zeta & = i \, \pi^2 \wedge \kappa + \pi^1 \wedge \zeta - \ov{\pi^1} \wedge \zeta,  \\
d \ov{\kappa} & = \ov{\pi^1} \wedge \ov{\kappa} + \ov{\pi^2} \wedge \rho - \kappa \wedge \ov{\zeta}, \\
d \ov{\zeta} & = - i \, \ov{\pi^2} \wedge \ov{\kappa} + \ov{\pi^1} \wedge \ov{\zeta} - {\pi^1} \wedge \ov{\zeta}.  \\
\end{aligned}
\end{equation}
We could obtain 
the expressions of $d \pi^1$ 
and $d \pi^2$ by taking the exterior derivative of the previous five equations.
But for now, as we have explicit expressions of $\pi^1$ 
and $\pi^2$ given by formulae $(\ref{eq:deltatilde})$ and $(\ref{eq:pi})$, we can perform an actual computation:

\begin{multline*}
d \pi^1 =  d {\sf t} \wedge \rho \\
+ X^1_{\rho \kappa} \, 
\rho \wedge \kappa + X^1_{\rho \zeta} \, 
\rho \wedge \zeta + X^1_{\rho \ov{\kappa}} \, \rho \wedge \ov{\kappa} + X^1_{\rho \ov{ \zeta}} \, 
\rho \wedge \ov{\zeta} \\
+ X^1_{\rho \pi^1} \, \rho \wedge \pi^1  + X^1_{\rho \pi^2} \, \rho \wedge \pi^2 
+ X^1_{\rho \ov{\pi^1}} \, \rho \wedge \ov{\pi^1} \\ +  X^1_{\rho \ov{\pi^2}} \, \rho \wedge \ov{\pi^2} 
+ i  \, \kappa \wedge \ov{\pi^2} + \zeta \wedge \ov{\zeta}         
,\end{multline*}
and
\begin{multline*}
d \pi^2 =  d{\sf t} \wedge \kappa \\
+ X^2_{\rho \kappa} \, 
\rho \wedge \kappa + X^2_{\kappa \zeta} \,
 \kappa \wedge \zeta + X^2_{\kappa \ov{\kappa}} \, \kappa \wedge \ov{\kappa} + X^2_{\kappa \ov{ \zeta}} \, 
\kappa \wedge \ov{\zeta} \\
+ X^2_{\kappa \pi^1} \, \kappa \wedge \pi^1  + X^2_{\kappa \pi^2} \, \kappa \wedge \pi^2 
+ X^2_{\kappa \ov{\pi^1}} \, \kappa \wedge \ov{\pi^1} \\  +  X^2_{\kappa \ov{\pi^2}} \, \kappa \wedge \ov{\pi^2} 
+  \zeta  \wedge \ov{\pi^2} + \pi^2 \wedge \ov{\pi^1}         
.\end{multline*}

From these equations, we see 
that the absorption is straightforward and that there remain no nonconstant essential torsion term. Indeed if we define the 
absorbed form $\Lambda$ by:
\begin{equation*}
\Lambda = d {\sf t} - X^2_{\rho \kappa} \, \rho -
 X^1_{\rho \kappa} \,\kappa - \sum_{\nu = \zeta, \pi^1, \cdots, \ov{\pi^2}} X^1_{\rho \nu} \, \nu,
\end{equation*}
the previous two equations become:
\begin{equation*}
d \pi^1 =  \Lambda \wedge \rho + i  \, \kappa \wedge \ov{\pi^2} + \zeta \wedge \ov{\zeta},        
\end{equation*}
and
\begin{equation*}
d \pi^2 =  \Lambda \wedge \kappa
+  \zeta  \wedge \ov{\pi^2} + \pi^2 \wedge \ov{\pi^1}.        
\end{equation*}
A straightforward computation gives the expression of $d \Lambda$:
\begin{equation*}
d \Lambda = - \pi^1 \wedge \Lambda + i\, \pi^2 \wedge \ov{\pi^2} - \pi^1 \wedge \ov{\Lambda}.
\end{equation*}
Let us summarize the results that we have obtained so far:
The ten $1$-forms $\rho$, $\kappa$, $\zeta$, $\ov{\kappa}$,
$\ov{\zeta}$, $\pi^1$, $\pi^2$, $\ov{\pi^1}$, $\ov{\pi^2}$, $\Lambda$ satisfies the structure
equations given by $(\ref{eq:structure})$. This completes the proof ofh Theorem~\ref{thm:LC}.

\appendix
\begin{small}
\section{Torsion coefficients for the $G$-structures on ${\sf B }$}
\subsection{Coefficients $U^{\bullet}_{\bullet \bullet}$}

\begin{equation*}
U^{\sigma}_{\sigma \rho} = {\frac {\ee}{{\A}^{3}}}+{\frac {\dd}{{\A}^{3}}}
,\end{equation*}

\begin{equation*}
U^{\sigma}_{\sigma \zeta} = -{\frac {\cc}{{\A}^{3}}}
,\end{equation*}

\begin{equation*}
U^{\sigma}_{\sigma \ov{\zeta}} = -{\frac {\cc}{{\A}^{3}}}
,\end{equation*}

\begin{equation*}
U^{\rho}_{\sigma \rho} = {\frac {\cc \ee}{{\A}^{6}}}+{\frac {\cc\dd}{{\A}^{6}}} - {\frac {i\bb\ee}{{\A}^{5}}}+{\frac {i\dd{ \bbb}}
{{\A}^{5}}} 
,\end{equation*}

\begin{equation*}
U^{\rho}_{\sigma \zeta} = {\frac {i\ee}{{\A}^{3}}}-{\frac {i{ \bbb}\,\cc}{{\A}^{5}}}-{\frac {{\cc}^{2}}
{{\A}^{6}}}
,\end{equation*}

\begin{equation*}
U^{\rho}_{\sigma \ov{\zeta}} = {\frac {i\bb\cc}{{\A}^{5}}}-{\frac {{\cc}^{2}}{{\A}^{6}}}-{\frac {i\dd}{{\A}^{3}}
}
,\end{equation*}

\begin{equation*}
U^{\rho}_{\rho \zeta} = {\frac {\cc}{{\A}^{3}}}+{\frac {i{ \bbb}}{{\A}^{2}}}
,\end{equation*}

\begin{equation*}
U^{\rho}_{\rho \ov{\zeta}} = {\frac {\cc}{{\A}^{3}}}-{\frac {i\bb}{{\A}^{2}}} 
,\end{equation*}

\begin{equation*}
U^{\zeta}_{\sigma \rho} = {\frac {{\dd}^{2}}{{\A}^{6}}}+{\frac {i{ 
\bbb}\,\dd\bb}{{\A}^{7}}}-{\frac {i\ee{\bb}^{2}}{{\A}^{7}}}+{\frac {\dd\ee}{{\A}^{6}}}
,\end{equation*}

\begin{equation*}
U^{\zeta}_{\sigma \zeta} = {\frac {i\bb\ee}{{\A}^{5}}}-{\frac {i{ \bbb}\,\cc\bb}{{\A}^{7}}}-{\frac {\cc\dd}{{\A
}^{6}}}
,\end{equation*}

\begin{equation*}
U^{\zeta}_{\sigma \ov{\zeta}} = -{\frac {\cc\dd}{{\A}^{6}}}-{\frac {i\dd\bb}{{\A}^{5}}}+{
\frac {i{\bb}^{2}\cc}{{\A}^{7}}} 
,\end{equation*}

\begin{equation*}
U^{\zeta}_{\rho \zeta}  = {\frac {\dd}{{\A}^{3}}}+{\frac {i{ \bbb}\,\bb}{{\A}^{4}}} 
,\end{equation*}

\begin{equation*}
U^{\zeta}_{\rho \ov{\zeta}}= {\frac {\dd}{{\A}^{3}}}-{\frac {i{\bb}^{2}}{{\A}^{4}}}
,\end{equation*}

\begin{equation*}
U^{\zeta}_{\zeta \overline{\zeta}} = {\frac {i\bb}{{\A}^{2}}}
.\end{equation*}

\newpage

\section{Torsion coefficients for the $G$-structures on ${\sf N}$}
\subsection{
Coefficients $U^{\bullet}_{\bullet \bullet}$}

\begin{equation*}
U^{\tau}_{\tau \sigma} = \frac{\hh}{\A^4} - \frac{\bbb \G}{\A^6} - \frac{\bb \G}{\A^6} + \frac{\kk}{\A^6}
,\end{equation*}

\begin{equation*}
U^{\tau}_{\tau \rho} =  {\frac {\bb\ff}{{\A}^{6}}}+{\frac {{ \bbb}\,\ff}{{\A}^{6}}}
,\end{equation*}

\begin{equation*}
U^{\tau}_{\tau \zeta} = -{\frac {\ff}{{\A}^{4}}}
,\end{equation*}

\begin{equation*}
U^{\tau}_{\tau \ov{\zeta}} =  -{\frac {\ff}{{\A}^{4}}}
,\end{equation*}

\begin{equation*}
U^{\tau}_{\sigma \rho} =  -{\frac {\bb}{{\A}^{2}}}-{\frac {{ \bbb}}{{\A}^{2}}}
,\end{equation*}

\begin{equation*}
U^{\sigma}_{\tau \sigma} = {\frac {\G\ee}{{\A}^{7}}}-{\frac {\hh\cc}{{\A}^{7}}}-{\frac {\kk\cc}{{\A}^{7}}}+{
\frac {\G\dd}{{\A}^{7}}}+{\frac {\ff\kk}{{\A}^{8}}}+{\frac {\ff\hh}{{\A}^{8}}}-{
\frac {\ff{ \bbb}\,\G}{{\A}^{10}}}-{\frac {\ff\bb\G}{{\A}^{10}}}
,\end{equation*}

\begin{equation*}
U^{\sigma}_{\tau \rho} = 
{\frac {{ \bbb}\,{\ff}^{2}}{{\A}^{10}}}+{\frac {\bb{\ff}^{2}}{{\A}^{10}}}-{
\frac {\ff\ee}{{\A}^{7}}}-{\frac {\ff\dd}{{\A}^{7}}}+{\frac {\kk}{{\A}^{4}}}+{
\frac {\hh}{{\A}^{4}}}
,\end{equation*}

\begin{equation*}
U^{\sigma}_{\tau \zeta} = -{\frac {\G}{{\A}^{4}}}+{\frac {\cc\ff}{{\A}^{7}}}-{\frac {{\ff}^{2}}{{\A}^{8}}} 
,\end{equation*}

\begin{equation*}
U^{\sigma}_{\tau \ov{\zeta}}  = -{\frac {\G}{{\A}^{4}}}+{\frac {\cc\ff}{{\A}^{7}}}-{\frac {{\ff}^{2}}{{\A}^{8}}}
,\end{equation*}

\begin{equation*}
U^{\sigma}_{\sigma \rho} = {\frac {\ee}{{\A}^{3}}}+{\frac {\dd}{{\A}^{3}}}-{\frac {\bb\ff}{{\A}^{6}}}-{
\frac {{ \bbb}\,\ff}{{\A}^{6}}}
,\end{equation*}

\begin{equation*}
U^{\sigma}_{\sigma \zeta} = -{\frac {\cc}{{\A}^{3}}}+{\frac {\ff}{{\A}^{4}}}
,\end{equation*}

\begin{equation*}
U^{\sigma}_{\sigma \ov{\zeta}} = -{\frac {\cc}{{\A}^{3}}}+{\frac {\ff}{{\A}^{4}}}
,\end{equation*}

\begin{equation*}
U^{\rho}_{\tau \sigma} = {\frac {-i\ee\bb\G}{{\A}^{9}}}-{\frac {i{ \bbb}\,\cc\hh}{{\A}^{9}}}+{\frac {i\dd{
 \bbb}\,\G}{{\A}^{9}}}+{\frac {i\bb\cc\kk}{{\A}^{9}}}+{\frac {\ee\G\cc}{{\A}^{10}}}+
{\frac {\dd\G\cc}{{\A}^{10}}}-{\frac {i \dd \kk}{{\A}^{7}}}+{\frac {i\ee\hh}{{\A}^{7}}}-
{\frac {{\cc}^{2}\hh}{{\A}^{10}}}-{\frac {{\cc}^{2}\kk}{{\A}^{10}}}-{\frac {{\G}^
{2}\bb}{{\A}^{10}}}+{\frac {\G\kk}{{\A}^{8}}}+{\frac {\G\hh}{{\A}^{8}}}-{\frac {{
\G}^{2}{ \bbb}}{{\A}^{10}}} 
,\end{equation*}

\begin{equation*}
U^{\rho}_{\tau \rho} = -{\frac {\cc\dd\ff}{{\A}^{10}}}+{\frac {\ff\bb\G}{{\A}^{10}}}+{\frac {\ff{ \bbb}\,\G}
{{\A}^{10}}}-{\frac {\cc\ee\ff}{{\A}^{10}}}+{\frac {\hh\cc}{{\A}^{7}}}+{\frac {\kk\cc}{
{\A}^{7}}}-{\frac {i\bb\kk}{{\A}^{6}}}+{\frac {i\bb\ee\ff}{{\A}^{9}}}-{\frac {i\dd{
 \bbb}\,\ff}{{\A}^{9}}}+{\frac {i{ \bbb}\,\hh}{{\A}^{6}}} 
,\end{equation*}

\begin{equation*}
U^{\rho}_{\tau \zeta} = {\frac {i{ \bbb}\,\cc\ff}{{\A}^{9}}}-{\frac {i\ee\ff}{{\A}^{7}}}-{\frac {i{ 
\bbb}\,\G}{{\A}^{6}}}+{\frac {{\cc}^{2}\ff}{{\A}^{10}}}-{\frac {\G\cc}{{\A}^{7}}}+{
\frac {i\kk}{{\A}^{4}}}-{\frac {\G\ff}{{\A}^{8}}} 
,\end{equation*}

\begin{equation*}
U^{\rho}_{\tau \overline{\zeta}} = {\frac {-i\bb\cc\ff}{{\A}^{9}}}+{\frac {i\dd\ff}{{\A}^{7}}}+{\frac {i\bb\G}{{\A}^{6}}}
+{\frac {{\cc}^{2}\ff}{{\A}^{10}}}-{\frac {\G\cc}{{\A}^{7}}}-{\frac {\G\ff}{{\A}^{8
}}}-{\frac {i\hh}{{\A}^{4}}}
,\end{equation*}

\begin{equation*}
U^{\rho}_{\sigma \rho} = {\frac {\cc \ee}{{\A}^{6}}}+{\frac {\cc\dd}{{\A}^{6}}}-{\frac {\G\bb}{{\A}^{6}}}-{
\frac {\G{ \bbb}}{{\A}^{6}}}-{\frac {i\bb\ee}{{\A}^{5}}}+{\frac {i\dd{ \bbb}}
{{\A}^{5}}} 
,\end{equation*}

\begin{equation*}
U^{\rho}_{\sigma \zeta} = {\frac {i\ee}{{\A}^{3}}}-{\frac {i{ \bbb}\,\cc}{{\A}^{5}}}-{\frac {{\cc}^{2}}
{{\A}^{6}}}+{\frac {\G}{{\A}^{4}}} 
,\end{equation*}

\begin{equation*}
U^{\rho}_{\sigma \ov{\zeta}} = {\frac {i\bb\cc}{{\A}^{5}}}-{\frac {{\cc}^{2}}{{\A}^{6}}}-{\frac {i\dd}{{\A}^{3}}
}+{\frac {\G}{{\A}^{4}}} 
,\end{equation*}

\begin{equation*}
U^{\rho}_{\rho \zeta} = {\frac {\cc}{{\A}^{3}}}+{\frac {i{ \bbb}}{{\A}^{2}}}
,\end{equation*}

\begin{equation*}
U^{\rho}_{\rho \ov{\zeta}} = {\frac {\cc}{{\A}^{3}}}-{\frac {i\bb}{{\A}^{2}}} 
,\end{equation*}

\begin{dmath*}
U^{\zeta}_{\tau \sigma} = {\frac {-i\ee{\bb}^{2}\G}{{\A}^{11}}}-{\frac {i\dd\kk\bb}{{\A}^{9}}}+{\frac {i\ee\hh\bb}{
{\A}^{9}}}+{\frac {i{\bb}^{2}\cc\kk}{{\A}^{11}}}+{\frac {\hh\kk}{{\A}^{8}}}+{\frac 
{{\dd}^{2}\G}{{\A}^{10}}}-{\frac {i{ \bbb}\,\cc\hh\bb}{{\A}^{11}}}+{\frac {i\dd{
 \bbb}\,\G\bb}{{\A}^{11}}}-{\frac {\cc\dd\kk}{{\A}^{10}}}-{\frac {\cc\dd\hh}{{\A}^{10}}
}+{\frac {\G\ee\dd}{{\A}^{10}}}-{\frac {\hh{ \bbb}\,\G}{{\A}^{10}}}-{\frac {\hh\bb\G
}{{\A}^{10}}}+{\frac {{\hh}^{2}}{{\A}^{8}}}
,\end{dmath*}

\begin{equation*}
U^{\zeta}_{\tau \rho} = {\frac {\kk\dd}{{\A}^{7}}}-{\frac {\dd\ee\ff}{{\A}^{10}}}+{\frac {\hh{ \bbb}\,\ff}{{\A
}^{10}}}+{\frac {\hh\bb\ff}{{\A}^{10}}}-{\frac {i\kk{\bb}^{2}}{{\A}^{8}}}-{\frac {
i{ \bbb}\,\dd\ff\bb}{{\A}^{11}}}+{\frac {i{ \bbb}\,\hh\bb}{{\A}^{8}}}+{\frac {i\ee
\ff{\bb}^{2}}{{\A}^{11}}}+{\frac {\dd\hh}{{\A}^{7}}}-{\frac {{\dd}^{2}\ff}{{\A}^{10}}
} 
,\end{equation*}

\begin{equation*}
U^{\zeta}_{\tau \zeta} = -{\frac {\G\dd}{{\A}^{7}}}-{\frac {\ff\hh}{{\A}^{8}}}-{\frac {i\bb\ee\ff}{{\A}^{9}}}-{
\frac {i{ \bbb}\,\G\bb}{{\A}^{8}}}+{\frac {i{ \bbb}\,\cc\ff\bb}{{\A}^{11}}}+{
\frac {\cc\dd\ff}{{\A}^{10}}}+{\frac {i\bb\kk}{{\A}^{6}}} 
,\end{equation*}

\begin{equation*}
 U^{\zeta}_{\tau \ov{\zeta}}= -{\frac {\G\dd}{{\A}^{7}}}-{\frac {\ff\hh}{{\A}^{8}}}+{\frac {i\dd\ff\bb}{{\A}^{9}}}+{
\frac {\cc\dd\ff}{{\A}^{10}}}-{\frac {i\hh\bb}{{\A}^{6}}}+{\frac {i{\bb}^{2}\G}{{\A}^{
8}}}-{\frac {i{\bb}^{2}\cc\ff}{{\A}^{11}}} 
,\end{equation*}

\begin{equation*}
U^{\zeta}_{\sigma \rho} = {\frac {{\dd}^{2}}{{\A}^{6}}}-{\frac {\hh{ \bbb}}{{\A}^{6}}}+{\frac {i{ 
\bbb}\,\dd\bb}{{\A}^{7}}}-{\frac {i\ee{\bb}^{2}}{{\A}^{7}}}+{\frac {\dd\ee}{{\A}^{6}}}-
{\frac {\hh\bb}{{\A}^{6}}} 
,\end{equation*}

\begin{equation*}
U^{\zeta}_{\sigma \zeta} = {\frac {i\bb\ee}{{\A}^{5}}}-{\frac {i{ \bbb}\,\cc\bb}{{\A}^{7}}}-{\frac {\cc\dd}{{\A
}^{6}}}+{\frac {\hh}{{\A}^{4}}}
,\end{equation*}

\begin{equation*}
U^{\zeta}_{\sigma \ov{\zeta}} = -{\frac {\cc\dd}{{\A}^{6}}}-{\frac {i\dd\bb}{{\A}^{5}}}+{\frac {\hh}{{\A}^{4}}}+{
\frac {i{\bb}^{2}\cc}{{\A}^{7}}} 
,\end{equation*}

\begin{equation*}
U^{\zeta}_{\rho \zeta}  = {\frac {\dd}{{\A}^{3}}}+{\frac {i{ \bbb}\,\bb}{{\A}^{4}}} 
,\end{equation*}

\begin{equation*}
U^{\zeta}_{\rho \ov{\zeta}}= {\frac {\dd}{{\A}^{3}}}-{\frac {i{\bb}^{2}}{{\A}^{4}}}
,\end{equation*}

\begin{equation*}
U^{\zeta}_{\zeta \overline{\zeta}} = {\frac {i\bb}{{\A}^{2}}}.
\end{equation*}

\section{Torsion coefficients for the $G$-structures on $\LC$}
\subsection{Coefficients $T^{\bullet}_{\bullet \bullet}$} 

\begin{equation*}
T^{\rho}_{\rho \kappa} = i\, {\frac { \overline{\sf b}}{{\sf c}\overline{ \sf c}}} - {\frac {\sf e  }{{\sf c}{ \sf f}}}\, 
 \frac{\ov{z_2}}{1 - z_2 \ov{z_2}},
\end{equation*}

\begin{equation*}
T^{\rho}_{\rho \zeta} =  \frac{1}{\ff} \, \frac{\ov{z_2}}{1 - z_2 \ov{z_2} },
\end{equation*}

\begin{equation*}
T^{\rho}_{\rho \overline{\kappa}}
=
-i \, {\frac {\sf b}{{\sf c} \overline{ \sf c}}}
+
{\frac {\overline{\sf e}}{\overline{ \sf c}\overline{\sf f}}} \, - \frac{z_2}{1-z_2 \ov{z_2}},
\end{equation*}

\begin{equation*}
T^{\rho}_{\rho \overline{\zeta}} = \frac {1}{\fb} \,  \frac{z_2}{1-z_2 \ov{z_2}}
,
\end{equation*}

\begin{equation*}
T^{\kappa}_{\rho \kappa}
=
\frac {{\sf e}
\overline{\sf b} }{{\sf c} \overline{ \sf c}^2{\sf f}} \, \frac{1}{1-z_2 \ov{z_2}}
+
\frac {{\sf d} }{{\sf c}\overline{ \sf c}{\sf f}}\,  \frac{\ov{z_2}}{1 - z_2 \ov{z_2}} +
 i\, {\frac
{{\sf b}\overline{\sf b}}{{\sf c}^{2} \overline{ \sf c} ^{2}}}
-
{\frac {{\sf e} {\sf b}}{{\sf c}^{2}\overline{ \sf c}{\sf f}}}\,  \frac{\ov{z_2}}{1 - z_2 \ov{z_2}}
\end{equation*}

\begin{equation*}
T^{\kappa}_{\rho \zeta}
=
-{\frac {\overline{\sf b}}{
\overline{ \sf c} ^{2}{\sf f} }}\,\frac{1}{1-z_2 \ov{z_2}} 
,\end{equation*}

\begin{equation*}
T^{\kappa}_{\rho \overline{\kappa}}
=
\frac {\sf d }{\overline{ \sf c}
  ^2 \ff}\,\frac{1}{1-z_2 \ov{z_2}}
-
 \frac {{\sf e} {\sf b}}{{\sf c} 
\overline{ {\sf c}} ^2 \ff }\,\frac{1}{1-z_2 \ov{z_2}} 
-
i\, {\frac {{\sf b}^{2}}{{\sf c}^{2} 
\overline{ \sf c} ^{2}}}
+
{\frac {{\sf b}\overline{\sf e} }{{\sf c} \overline{ \sf c} ^{2}\overline{\sf f}}\, - \frac{z_2}{1-z_2 \ov{z_2}}
}
,\end{equation*}

\begin{equation*}
T^{\kappa}_{\rho \overline{\zeta}}
=
\frac {  \sf b}{{\sf c}\overline{ \sf c}
\overline{\sf f}} \, \frac{z_2}{1-z_2 \ov{z_2}},
 \end{equation*}

\begin{equation*}
T^{\kappa}_{\kappa \zeta}
=
\frac{1}{\ff} \,   \frac{\ov{z_2}}{1 - z_2 \ov{z_2}}
,
\end{equation*}

\begin{equation*}
T^{\kappa}_{\kappa \overline{\kappa}}
=
{\frac {\sf e}{\overline{ \sf c} \ff}}\,\frac{1}{1-z_2 \ov{z_2}} 
+
i\, {\frac {\sf b}{
{\sf c}\overline{ \sf c}}}
,\end{equation*}

\begin{equation*}
T^{\kappa}_{\zeta \overline{\kappa}}
=
-\frac {{\sf c}}{\overline{ \sf c} \ff}\,\frac{1}{1-z_2 \ov{z_2}} 
,\end{equation*}

\begin{equation*}
T^{\zeta}_{\rho \kappa}
=
\frac {{\sf 
e}^{2}\overline{\sf b}}{{\sf c}^{2}
\overline{ \sf c}  ^2 \ff} \, \frac{1}{1-z_2 \ov{z_2}} 
+
i\, {\frac {{\sf d}\overline{\sf b}}{{\sf c}^{2} 
\overline{ \sf c} ^2}}
,\end{equation*}

\begin{equation*}
T^{\zeta}_{\rho \zeta}
=
-
{\frac { {\sf e}\overline{\sf b}}{{\sf c} \overline{ \sf c} ^{2} \ff}}\, \frac{1}{1-z_2 \ov{z_2}} 
-
{\frac { {\sf e} {\sf b}}{{\sf c}^{2}
\overline{ \sf c}{\sf f}}} \,  \frac{\ov{z_2}}{1 - z_2 \ov{z_2}}
+
{\frac {{ \sf d} }{{\sf c}\overline{ \sf c}{\sf f}}}\,  \frac{\ov{z_2}}{1 - z_2 \ov{z_2}}
,
\end{equation*}

\begin{equation*}
T^{\zeta}_{\rho \overline{\kappa}}
=
\frac { {\sf e}{ \sf d}}{\cc \cb^2 \ff }\, \frac{1}{1-z_2 \ov{z_2}}
-
{\frac {{\sf e}^{2}{\sf b}}{{\sf c}^{2
}\overline{ \sf c}^{2} \ff}} \, \frac{1}{1-z_2 \ov{z_2}}
-
i\,{\frac {{\sf b}{\sf d}}{{\sf c}^{2} 
\overline{ \sf c}  ^{2}}}
-
{\frac {{ \sf d}\overline{\sf e}}{{\sf c}  \overline{ \sf c}^{2}\overline{\sf f}} \,  \frac{z_2}{1-z_2 \ov{z_2}}
}
,\end{equation*}

\begin{equation*}
T^{\zeta}_{\rho \overline{\zeta}}
=
{\frac {{ \sf d}}{{\sf c}\overline{ \sf c}
\overline{\sf f}}}\, \frac{z_2}{1 - z_2 \ov{z_2}}
,\end{equation*}

\begin{equation*}
T^{\zeta}_{\kappa \zeta}
=
+{\frac {{ \sf e}}{{\sf c}\overline{ \sf c}
\overline{\sf f}}}\, \frac{\ov{z_2}}{1 - z_2 \ov{z_2}} 
,
\end{equation*}

\begin{equation*}
T^{\zeta}_{\kappa \overline{\kappa}}
=
{\frac {{\sf e}^{2}}{{\sf c}\overline{ \sf c}{\sf f}}}\,\frac{1}{1-z_2 \ov{z_2}}
+
{
i\,\frac {\sf d}{{\sf c}\overline{ \sf c}}}
,\end{equation*}

\begin{equation*}
T^{\zeta}_{\zeta \overline{\kappa}}
=
-\frac { \ee}{ \cb \ff} \, \frac{1}{1-z_2 \ov{z_2}}
.\end{equation*}

\subsection{Coefficients $U^{\bullet}_{\bullet \bullet}$}
\[
U^{\rho}_{\rho \kappa} = i\,  {\frac {\overline{\sf b}}{{\sf c}\overline{ \sf c}}}
+
{\frac { {\sf e}\overline{ \sf c}}{{\sf c}^{2}}}\,\ov{z_2}
,\]

\[
U^{\rho}_{\rho \zeta} = -
{\frac {\overline{ \sf c}}{{\sf c}}}\, \ov{z_2}
,\]

\[
U^{\rho}_{\rho \overline{\kappa}}
=
-i {\frac {{\sf b}}{{\sf c}\overline{ \sf c}}}
+
{\frac {\overline{\sf e}{\sf c}}{  \overline{ \sf c}^2}} \, z_2
,\]

\[
U^{\rho}_{\rho \overline{\zeta}} = 
-
\frac{{\sf c}}{\overline{ \sf c}} \, z_2
,\]

\[
U^{\kappa}_{\rho \kappa}
=
-
{\frac { {\sf e}
\overline{\sf b}}{{\sf c}^{2}\overline{ \sf c}}}
-
{\frac {{ \sf d}}{{\sf c}^
{2}}} \, \ov{z_2}
+
i\, {\frac {{\sf b}\overline{\sf b}}{{\sf c}^{2}
  \overline{ \sf c}^2}}
  +
  {\frac {{\sf b} {\sf e}}{{\sf c}
^{3}}}\,\ov{z_2}
,\]

\[
U^{\kappa}_{\rho \zeta}
= \frac{\bbb}{\cc \cb}
,\]

\[
U^{\kappa}_{\rho \overline{\kappa}}
=
-{\frac {\sf d}{{\sf c}\overline{ \sf c}}}
+
{\frac { {\sf e}{\sf b}}{{\sf c}^{2}\overline{ \sf c}}}
-
i\, {\frac {
{\sf b}^{2}}{{\sf c}^{2}  \overline{ \sf c}^2}}
+
{\frac {{\sf b}
\overline{\sf e}}{  
\overline{ \sf c} ^{3}}}\, z_2
,\]

\[
U^{\kappa}_{\rho \overline{\zeta}}
=
-{\frac {\sf b}{ 
\overline{ \sf c}^{2}}}\, z_2
,\]

\[
U^{\kappa}_{\kappa \zeta}
=
-{\frac {\overline{ \sf c}}{{\sf c}}}\,\ov{z_2}
,\]

\[
U^{\kappa}_{\kappa \overline{\kappa}}
=
-{\frac {\sf e}{\sf c}}+i\, {\frac {{\sf b}}{{\sf c}\overline{ \sf c}}}
,\]

\begin{equation*}
U^{\zeta}_{\rho \kappa}
=
-
{\frac { {\sf e}\overline{\sf d}
}{{\sf c} \overline{ \sf c}
 ^{2}}}\, z_2
 +
 {\frac { \bbb {\sf e}\overline{\sf e}} { 
\overline{ \sf c}  ^{3}{\sf c}}} \, z_2
-
{\frac {
{\sf e}^{2}\overline{\sf b}}{\overline{ \sf c}{\sf c}^{3}}}
+
i\, {\frac {{ \sf d}\overline{\sf b}}{{\sf c}
^{2}  \overline{ \sf c}^2}}
,\end{equation*}

\begin{equation*}
U^{\zeta}_{\rho \zeta}
=
{\frac {\overline{\sf d}}{
  \overline{ \sf c}^2}} \, z_2
  -
  {
\frac {\overline{\sf e} \overline{\sf b}}{  \overline{ \sf c}^{3}}}\, z_2
+
\frac{\ee \bbb}{\cc^2 \cb}
+
{\frac {{\sf b}{\sf e}}{{\sf c}
^{3}}}\, \ov{z_2}
-
{\frac {{ \sf d}}{{\sf c}
^{2}}}\, \ov{z_2}
,\end{equation*}

\begin{equation*}
U^{\zeta}_{\rho \overline{\kappa}}
=
2\,{\frac {\overline{\sf e}{ \sf d}}{
  \overline{ \sf c}^{3}}} \, z_2
  -
  {\frac { {\sf e}\overline{\sf e} {\sf b}}{
  \overline{ \sf c}^{3}{\sf c}}} \, z_2
-
{\frac { {\sf e}{\sf d}}{{\sf c}^{2}\overline{ \sf c}}}
+ 
{\frac {{\sf e}^{2}{\sf b}}{\overline{ \sf c}{\sf c}
^{3}}}
-
i\, {\frac {{\sf d}{\sf b}}{{\sf c}^{2}  \overline{ \sf c}^2}}
,\end{equation*}

\begin{equation*}
U^{\zeta}_{\rho \overline{\zeta}}
=
-2\,{\frac { {\sf d} }{ 
\overline{ \sf c}^{2}}} \, z_2
+
{\frac {{\sf e} {\sf b}
}{{\sf c } \overline{ \sf c}^{2}}} \, z_2
,
\end{equation*}

\[
U^{\zeta}_{\kappa \zeta}
=
-
{\frac { {\sf e}\overline{ \sf c}}{{\sf c}^{2}}} \, \ov{z_2}
,\]

\[
U^{\zeta}_{\kappa \overline{\kappa}}
=
{\frac { {\sf e}\overline{\sf e} }{ \overline{ \sf c}^2
 }} \, z_2
-
{\frac {{\sf e}^{2}}{{\sf c}^{2}
}}
+
i\, {\frac {\sf d}{{\sf c}\overline{ \sf c}}}
,\]

\[
U^{\zeta}_{\kappa \overline{\zeta}}
=
- {\frac { {\sf e}}{\overline{ \sf c}}}\, z_2
,\]

\[
U^{\zeta}_{\zeta \overline{\kappa}}
=
- {\frac {\overline{\sf e}{\sf c}}{
  \overline{ \sf c}^2}} \, z_2
+
{\frac {\sf e}{\sf c}}
,\]

\[
U^{\zeta}_{\zeta \overline{\zeta}}
=
{\frac {{\sf c}}{\overline{ \sf c}}}\,z_2
.\]

\subsection{
Coefficients $V^{\bullet}_{\bullet \bullet}$} 

\begin{equation*}
V^{\rho}_{\rho \kappa}
=
-
\frac {\overline {\sf e}} {\overline{\sf c}}
+
{\frac {{\sf e}  \overline{\sf c} } {{{\sf c}}^{2} 
}} \, \ov{z_2}
,\end{equation*}

\begin{equation*}
V^{\rho}_{\rho \zeta} 
=
-
{\frac {{\overline{\sf  c}}  }{{\sf c} } \, \ov{z_2}}
,\end{equation*}

\begin{equation*}
V^{\rho}_{\rho \overline{\kappa}} 
=
-{\frac {{\sf e}} {{\sf c}}}
+
 {\frac {{\overline{\sf e}} {\sf c} }{{{\overline{\sf c}
}}^{2}}} \, z_2
,\end{equation*}

\begin{equation*}
V^{\rho}_{\rho \overline{\zeta}}
=
-\frac{\sf c}{\overline{\sf c}} 
  \, z_2
,\end{equation*}

\begin{equation*}
V^{\kappa}_{\rho \kappa}
=
-
{\frac {{ \sf d}}{{\sf c}^
{2}}} \, \ov{z_2}
- i \, 
  {\frac {\ee^2 \cb}{{\sf c}
^{3}}}\,\ov{z_2}
\end{equation*}

\begin{equation*}
V^{\kappa}_{\rho \zeta} 
=
i \, \frac{\overline{\sf e}}{\overline{\sf c}}
,\end{equation*}

\begin{equation*}
V^{\kappa}_{\rho \overline{\kappa}} 
=
-
{\frac {{\sf d}}{{\sf 
c} {\cb}}}
-
i \, \, \frac{{\sf e} \overline{\sf e} }{\overline{\sf c}^{2}} \, z_2
,\end{equation*}

\begin{equation*}
V^{\kappa}_{\rho \overline{\zeta}}
=
i \,  \frac {{\sf e} }{\overline{\sf c}} \, z_2
,\end{equation*}

\begin{equation*}
V^{\kappa}_{\kappa \zeta} 
=
- \frac{\overline{\sf c}}{\sf c} \, \ov{z_2}
,\end{equation*}

\begin{equation*}
V^{\kappa}_{\kappa \overline{\kappa}}
=0
,\end{equation*}

\begin{equation*}
V^{\zeta}_{\rho \kappa}
=
- \frac{\db \ee }{\cc \cb^2} \, z_2
-
{\frac {{\sf d}{\overline{\sf e}}}{{\sf c}
{{\overline{\sf c}}}^{2}}}
+
i \, \frac {{\sf e} \overline{\sf e}^{2}} {\overline{\sf c}^{3}} \, z_2
-
i \, \frac{{\sf e}^{2} \overline{\sf e}} {{\sf c}^{2}{\overline{\sf c}}}
,\end{equation*}

\begin{equation*}
V^{\zeta}_{\rho \zeta} 
=
-\frac{\db}{\cb^2} \, z_2
+
i \, {\frac {{\sf e}\,{\overline{\sf e}}}{{\sf c}
\,{\overline{\sf c}}}}
-
i \, \frac{\cb \ee^2} {\cc^3} \, \ov{z_2}
 -
i \, \frac{\cc \eb^2}{\cb^3} \, z_2
-
\frac{\sf d}{{\sf c}^{2}} \, \ov{z_2}
,\end{equation*} 

\begin{equation*}
V^{\zeta}_{\rho \overline{\kappa}} 
=
2 \, \frac{{\sf d} \overline{\sf e}} {{\overline{\sf c}}^{3}} \, z_2
 + 
i \, \frac{\overline{\sf e} {\sf e}^{2}} {
\overline{\sf c}^{2}{\sf c}} \, z_2
 -
2 \, {\frac {{\sf e} {\sf d}}{{\overline{\sf c}} {{\sf c
}}^{2}}}
-
i \, {\frac {{{\sf e}}^{3}}{{{\sf c}}^{3}}}
,\end{equation*}

\begin{equation*}
V^{\zeta}_{\rho \overline{\zeta}} 
=
-2 \, {\frac {{\sf d}}{{{
\overline{\sf c}}}^{2}}} \, z_2
- 
i \, {\frac {{{\sf e}}^{2}
}{{\sf c}\,{\overline{\sf c}}}} \, z_2
,\end{equation*}

\begin{equation*}
V^{\zeta}_{\kappa \zeta} 
=
-
{\frac {{\sf e}\,{\overline{\sf c}}}{{{\sf c}}^{2}}} \, \ov{z_2}
,\end{equation*}

\begin{equation*}
V^{\zeta}_{\kappa \overline{\kappa}}
=
{\frac {{\sf e}{\overline{\sf e}}}
{{{\overline{\sf c}}}^{2}}} \, z_2
-
{\frac {{{\sf e}}^{2}}{{{
\sf c}}^{2}}} 
+
i \, {\frac {{\sf d}}{{\sf c}\,{\overline{\sf c}}}}
,\end{equation*}

\begin{equation*}
V^{\zeta}_{\kappa \overline{\zeta}} 
=
- \frac{\sf e}{{\overline{\sf c}}} \, z_2
,\end{equation*}

\begin{equation*}
V^{\zeta}_{\zeta \overline{\kappa}}
=
-\frac{ \overline{\sf e} {\sf c}} { \overline{\sf c}^{2} } \, z_2
+
{\frac {{\sf e}}{{\sf c}}}
,\end{equation*}

\begin{equation*}
V^{\zeta}_{\zeta \overline{\zeta}}
=
{\frac{\sf c}{\overline{\sf c}}} \, z_2
.\end{equation*}

\subsection{Coefficients $X^{\bullet}_{\bullet \bullet}$}

\begin{equation*}
X^1_{\rho \kappa}= - \frac{1}{2} \, {\sf t} \, \frac{\cb \ee}{\cc} \, \ov{z_2} 
- \frac{3}{8} \, i \, \frac{\ee^2 \eb}{\cc^3} \, \ov{z_2} + \frac{1}{2} \, {\sf t} \frac{\eb}{\cb}
+ \frac{1}{8} \, i \, \frac{\ee \eb^2}{\cc \cb^2}
 \, z_2 \ov{z_2} +   \frac{1}{8} \, i \, \frac{\ee^3 \cb^2}{\cc^5} \, \ov{z_2}^2 +
 \frac{1}{4} \, i \, \frac{\eb^2 \ee}{\cc \cb^2}
,\end{equation*}

\begin{equation*}
X^1_{\rho \zeta} = - \frac{1}{4} \, i \, \frac{\eb^2}{\cb^2} \, + \frac{1}{2} \, i \, \frac{\ee \eb}{\cc^2} \, \ov{z_2} -
 \frac{1}{4} \, i \, \frac{\cb^2 \ee^2}{\cc^4} \, \ov{z_2}^2
,\end{equation*}

\begin{equation*}
X^1_{\rho \ov{\kappa}} = \ov{X^1_{\rho \kappa}}
,\end{equation*}

\begin{equation*}
X^1_{\rho \ov{\zeta}} = \ov{X^1_{\rho \zeta}}
,\end{equation*}

\begin{equation*}
X^1_{\rho \pi^1} = - {\sf t}
,\end{equation*}

\begin{equation*}
X^1_{\rho \pi^2} = \frac{1}{2} \, \frac{\eb}{\cb} + \frac{1}{2} \, \frac{\ee \cb}{\cc^2} \, \ov{z_2} 
,\end{equation*}

\begin{equation*}
X^1_{\rho \ov{\pi^1}} =  - {\sf t}
,\end{equation*}

\begin{equation*}
X^1_{\rho \ov{\pi^2}} = \ov{X^1_{\rho \pi^2}}
,\end{equation*}

\begin{equation*}
X^2_{\rho \kappa} = - \frac{1}{4} \, \frac{\ee \eb^3}{\cb^4} \, z_2 - \frac{1}{4} \frac{\ee^3 \eb}{\cb^4} \, \ov{z_2}
+ \frac{1}{8} \,  \, \frac{\ee^2 \eb^2}{\cc^2 \cb^2} 
\, \ov{z_2} + \frac{1}{4} \, \frac{\ee^2 \eb^2}{\cc^2 \cb^2} + \frac{1}{16} \, \frac{\eb^4 \cc^2}{\cb^6} \, z_2^2
+ \frac{1}{16} \, \frac{\ee^4 \cb^2}{\cc^6} \, \ov{z_2}^2 - {\sf t}^2, 
\end{equation*}

\begin{equation*}
X^2_{\kappa \nu } = X^1_{\rho \nu} \, \, \, \, \, \, \, \, \, \, \,  \text{ for} \, \, \,  \nu= \zeta, \pi^1, \cdots, \ov{\pi^2}.
\end{equation*}

\end{small}


\begin{thebibliography}{99}



{\bf\bibitem{BES-2007}
{\sf Beloshapka}}, V.K.; {\sf Ezhov}, V.; {\sf Schmalz}, G.:
{\em Canonical Cartan connection and holomorphic
invariants on Engel CR manifolds},
Russian J. Mathematical Physics {\bf 14} (2007), no.~2, 121--133.


{\bf\bibitem{EMS}
{\sf Ezhov}}, V.; {\sf McLaughlin}, B.; {\sf Schmalz}, G.:
{\em From Cartan to Tanaka: Getting Real in the Complex World},
Notices of the AMS (2011) {\bf 58}, no.~1.

{\bf\bibitem{Fels-Kaup-2007}
{\sf Fels}}, G.; {\sf Kaup}, W.:
{\em Classification of Levi-degenerate homogeneous CR-manifolds in dimension 
$5$}, Acta. Math. {\bf 201} (2008), 1--82. 

{\bf\bibitem{Isaev-Zaitsev}
{\sf Isaev}}, A.; {\sf Zaitsev}, D.:
{\em Reduction of five-dimensional uniformaly degenerate Levi CR structures to 
absolute parallelisms}, J. Anal. {\bf 23} (2013), no. 3, 1571--1605.

{\bf\bibitem{Kaup-Zaitsev}
{\sf Kaup}}, W.; {\sf Zaitsev}, D.: 
{\em On local CR-transformations of Levi-degenerate
group orbits in compact Hermitian symmetric spaces}, J. Eur. Math. Soc. {\bf 8} (2006),
465--490

{\bf\bibitem{Kobayashi}
{\sf Kobayashi}}, S.:
{\em Transformation groups in differential geometry}, 
Ergebnisse der Mathematik und ihrer Grenzgebiete, {\bf 70}, 
Springer-Verlag, New-York Heidelberg Berlin, 1972.

{\bf\bibitem{Medori-Spiro}
{\sf Medori}}, C.; {\sf Spiro}, A.:
{\em The equivalence problem for $5$-dimensional Levi degenerate 
CR manifolds}, Int. Math. Res. Not. (to appear), DOI: 10.1093/imrn/rnt129.

{\bf\bibitem{Merker-2003}
{\sf Gaussier}}, H.; {\sf Merker}, J.:
{\em A new example of uniformly Levi nondegenerate hypersurface
in $\C^3$}, Ark. Mat. {\bf 41} (2003), no.~1, 85--94.

{\bf\bibitem{MPS}
{\sf Merker}}, J.; {\sf Pocchiola}, S.; {\sf Sabzevari}, M.:
{\em Equivalences of $5$-dimensional CR-manifolds, II: General classes I, II, 
III-1, III-2, IV-1, IV-2}, arxiv.org/abs/1311.5669.

{\bf\bibitem{Olver-1995}
{\sf Olver}},~P.J.:
{\em Equivalence, Invariance and Symmetries}. Cambridge
University Press, Cambridge, 1995, xvi+525~pp.

{\bf\bibitem{pocchiola}
{\sf Pocchiola}}, S.:
{\em Absolute parallelism for $2$-nondegenerate real hypersurfaces $M^5 \subset \C^3$ of constant Levi rank $1$},
arxiv.org/abs/1312.6400, 56~pp. 


{\bf\bibitem{pocchiola2}
{\sf Pocchiola}}, S.:
{\em Canonical Cartan connection for $4$-dimensional CR-manifolds belonging to general class ${\sf II}$},
17~pp.

{\bf\bibitem{pocchiola3}
{\sf Pocchiola}}, S.:
{\em Canonical Cartan connection for $5$-dimensional CR-manifolds belonging to general class ${\sf III}_2$},
19~pp.


{\bf\bibitem{Sternberg}
{\sf Sternberg}}, S.:
{Lectures on Differential Geometry}. Prentice-Hall mathematical series.


\end{thebibliography}
\end{document}